\documentclass[aos,nameyear,seceqn,dvips]{arximspdf}

\doi{10.1214/07-AOS552}
\volume{37}
\issue{1}
\pubyear{2009}
\firstpage{73}
\lastpage{104}

\makeatletter

\def\varDelta{\Delta}
\def\varOmega{\Omega}
\def\varLambda{\Lambda}
\def\varXi{\Xi}

\newtheorem{lemma}{Lemma}
\newtheorem{proposition}{Proposition}

\newproclaim{remark}{Remark}
\newproclaim{example}{Example}

\newcommand{\fr}[1]{(\ref{#1})}
\newcommand{\ro}{\varrho}
\newcommand{\ph}{\varphi}
\newcommand{\af}{\alpha}
\newcommand{\eps}{\varepsilon}
\newcommand{\ga}{\gamma}
\newcommand{\te}{\theta}
\newcommand{\sgm}{\sigma}
\newcommand{\Lam}{\Lambda}
\newcommand{\hf}{\hat{f}_n}
\newcommand{\psimjk}{\psi_{mjk}}
\newcommand{\phmlk}{\ph_{mj_0k}}
\newcommand{\psijk}{\psi_{jk}}
\newcommand{\psiji}{\psi_{ji}}
\newcommand{\fm}{f_m}
\newcommand{\ym}{y_{m}}
\newcommand{\gm}{g_{m}}
\newcommand{\hmu}{h_{m}(u)}
\newcommand{\ymu}{y_{m}(u)}
\newcommand{\gma}{g_{m}(a)}
\newcommand{\bgma}{\overline{g_{m}(a)}}
\newcommand{\yma}{y_{m}(a)}
\newcommand{\gmu}{g_{m}(u)}
\newcommand{\bgmu}{\overline{g_{m}(u)}}
\newcommand{\zmu}{z_{m}(u)}
\newcommand{\ajk}{a_{j_0k}}
\newcommand{\bjk}{b_{jk}}
\newcommand{\hajk}{\widehat{a}_{j_0k}}
\newcommand{\hbjk}{\widehat{b}_{jk}}
\newcommand{\hfm}{\widehat{f}_m }
\newcommand{\Var}{\mbox{Var}}
\newcommand{\Bpqsa}{B_{p,q}^s (A) }
\newcommand{\Bpqs}{B_{p,q}^s}
\newcommand{\iab}{\int_a^b}
\newcommand{\intT}{\int_T}
\newcommand{\Aj}{A_j}
\newcommand{\Ujr}{U_{jr}}
\newcommand{\Bjr}{B_{jr}}
\newcommand{\hBjr}{\widehat{B}_{jr}}
\newcommand{\sumr}{\sum_{r \in\Aj}}
\newcommand{\sumk}{\sum_{k=0}^{2^j-1}}
\newcommand{\sumku}{\sum_{k \in\Ujr}}
\newcommand{\sumM}{\sum_{l=1}^M}
\newcommand{\sumN}{\sum_{i=1}^N}
\newcommand{\sumcj}{\sum_{m \in C_j}}

\newcommand{\sumjJ}{\sum_{j=j_0}^{J-1}}
\newcommand{\sumjone}{\sum_{j=j_0}^{j_1}}
\newcommand{\sumjtwo}{\sum_{j=j_1 +1}^{J-1}}
\newcommand{\sumjones}{\sum_{j=j_0}^{j_2}}
\newcommand{\sumjtwos}{\sum_{j=j_2 +1}^{J-1}}

\newcommand{\fjk}{f_{jk}}
\newcommand{\fjkp}{f_{jk'}}
\newcommand{\gaj}{\gamma_j}
\newcommand{\jo}{{j_0}}
\newcommand{\sstar}{{s^*}}
\newcommand{\EE}{\mathbb{E}}
\newcommand{\II}{\mathbb{I}}
\newcommand{\PP}{\mathbb{P}}

\newcommand{\ints}{\mathbb{Z}}
\newcommand{\RR}{\mathbb{R}}

\makeatother

\begin{document}
\begin{frontmatter}

\title{Functional deconvolution in a periodic setting: Uniform case}
\runtitle{Functional deconvolution}

\begin{aug}
\author[A]{\fnms{Marianna} \snm{Pensky}\thanksref{t1}\ead[label=e1]{mpensky@pegasus.cc.ucf.edu}}
and
\author[B]{\fnms{Theofanis} \snm{Sapatinas}\ead[label=e2]{T.Sapatinas@ucy.ac.cy}\corref{}}
\thankstext{t1}{Supported in part by NSF
Grants DMS-05-05133 and DMS-06-52524.}
\runauthor{M. Pensky and T. Sapatinas}
\affiliation{University of Central Florida and University of Cyprus}
\address[A]{Department of Mathematics \\
University of Central Florida \\
Orlando, Florida 32816-1353\\
USA \\
\printead{e1}} 
\address[B]{Department of Mathematics and Statistics \\
University of Cyprus \\
P.O. Box 20537 \\
Nicosia CY 1678\\
Cyprus \\
\printead{e2}}
\end{aug}

\received{\smonth{3} \syear{2007}}
\revised{\smonth{9} \syear{2007}}

%
\begin{abstract}
We extend deconvolution in a periodic setting to deal with
functional data. The resulting functional deconvolution model can be
viewed as a generalization of a multitude of inverse problems in
mathematical physics where one needs to recover initial or boundary
conditions on the basis of observations from a noisy solution of a
partial differential equation. In the case when it is observed at a
finite number of distinct points, the proposed functional
deconvolution model can also be viewed as a multichannel
deconvolution model.

We derive minimax lower bounds for the $L^2$-risk in the proposed
functional deconvolution model when $f(\cdot)$ is assumed to belong
to a Besov ball and the blurring function is assumed to possess some
smoothness properties, including both regular-smooth and
super-smooth convolutions. Furthermore, we propose an adaptive
wavelet estimator of $f(\cdot)$ that is asymptotically optimal (in
the minimax sense), or near-optimal within a logarithmic factor, in
a wide range of Besov balls.

In addition, we consider a discretization of the proposed functional
deconvolution model and investigate when the availability of
continuous data gives advantages over observations at the
asymptotically large number of points. As an illustration, we
discuss particular examples for both continuous and discrete
settings.
\end{abstract}

%


\pdfkeywords{62G05, 62G08, 35J05,
35K05,
35L05, Adaptivity,
Besov spaces,
block thresholding,
deconvolution,
Fourier analysis,
functional data,
Meyer wavelets,
minimax estimators,
multichannel deconvolution,
partial differen}

\end{frontmatter}

\section{Introduction}
\label{sec:intro}

We consider the estimation problem of the unknown response function
$f(\cdot)$ based on observations from the following noisy
convolutions:
%
\begin{equation}\label{conv1}
Y(u, t) = f*G (u, t) + \frac{\sigma(u)}{\sqrt{n}}
z(u,t), \qquad   u \in U,  \  t \in T,
\end{equation}
where
$U=[a,b]$, $-\infty< a \leq b < \infty$, and $T=[0,1]$. Here, $z(u,
t)$ is assumed to be a two-dimensional Gaussian white noise, that is, a
generalized two-dimensional Gaussian field with covariance function
\[
\EE[z (u_1, t_1) z (u_2, t_2)] = \delta(u_1-u_2) \delta(t_1-t_2),
\]
where $\delta(\cdot)$ denotes the Dirac $\delta$-function, $\sgm
(\cdot)$ is assumed to be a known positive function, and
%
\begin{equation}\label{conv2}
f*G (u,t) = \int_T f(x)G(u,t-x) \,dx,
\end{equation}
with the blurring (or
kernel) function $G(\cdot,\cdot)$ in \fr{conv2} also assumed to be
known. Note that, since $\sigma(\cdot)$ is assumed to be known,
both sides of \fr{conv1} can be divided by $\sigma(\cdot)$
leading to the equation
%
\begin{equation}
\label{convcont}
y(u, t) = \int_T f(x)g(u,t-x) \,dx +
\frac{1}{\sqrt{n}} z(u, t), \qquad  u \in U, \   t \in T,
\end{equation}
where $y(u, t)= Y(u, t)/\sigma(u)$ and
$g(u,t-x)= G(u,t-x)/\sigma(u)$. Consequently, without loss of
generality, we consider only the case when $\sigma(\cdot) \equiv1$
and thus, in what follows, we work with observations from model
(\ref{convcont}).

The model (\ref{convcont}) can be viewed as a \textit{functional
deconvolution} model. If $a=b$, it reduces to the standard
deconvolution model which attracted attention of a number of
researchers. After a rather rapid progress in this problem in late
1980s to early 1990s, authors turned to wavelet solutions of the
problem [see, e.g., Donoho (\citeyear{Don1995}), Abramovich and Silverman (\citeyear{AbrSil1998}),
Kalifa and Mallat (\citeyear{KalMal2003}), Johnstone, Kerkyacharian, Picard and
Raimondo (\citeyear{Joh2004}), Donoho and Raimondo (\citeyear{DonRai2004}), Johnstone and Raimondo
(\citeyear{JohRai2004}), Neelamani, Choi and Baraniuk (\citeyear{NeeChoBar2004}) and Kerkyacharian, Picard
and Raimondo (\citeyear{KerPicRai2007})]. The main effort was spent on producing adaptive
wavelet estimators that are asymptotically optimal (in the minimax
sense), or near-optimal within a logarithmic factor, in a wide range
of Besov balls and under mild conditions on the blurring function.
[For related results on the density deconvolution problem, we refer
to, e.g., Pensky and Vidakovic (\citeyear{PenVid1999}), Walter and Shen (\citeyear{WalShe1999}), Fan and
Koo (\citeyear{FanKoo2002}).]

On the other hand, the functional deconvolution model
(\ref{convcont}) can be viewed as a generalization of a multitude of
inverse problems in mathematical physics where one needs to recover
initial or boundary conditions on the basis of observations of a
noisy solution of a partial differential equation. Lattes and Lions
(\citeyear{LatLio1967}) initiated research in the problem of recovering the initial
condition for parabolic equations based on observations in a
fixed-time strip. This problem and the problem of recovering the
boundary condition for elliptic equations based on observations in
an internal domain were studied in Golubev and Khasminskii (\citeyear{GolKha1999});
the latter problem was also discussed in Golubev (\citeyear{Gol2004}). These and
other specific models are discussed in Section \ref{applic}.

Consider now a discretization of the functional deconvolution model
\fr{convcont} when $y(u,t)$ is observed at $n = N M$ points $(u_l,
t_i)$, $l=1,2,\ldots, M$, $i=1,2,\ldots, N$, that is,
%
\begin{equation}\label{convdis}
y(u_l, t_i) = \int_T f(x) g(u_l, t_i-x) \,dx + \eps_{li}, \qquad  u_l \in U, \   t_i
= i/N,
\end{equation}
where $\eps_{li}$ are standard Gaussian
random variables, independent for different $l$ and~$i$. In this
case, the functional deconvolution model \fr{convcont} can also be
viewed as a multichannel deconvolution problem considered in, for example,
Casey and Walnut (\citeyear{CasWal1994}) and De Canditiis and Pensky (\citeyear{DeCanPen2004}, \citeyear{DeCanPen2006}); this
model is also discussed in Section \ref{applic}.

Note that using the same $n$ in \fr{convcont} (continuous model) and
\fr{convdis} (discrete model) is not accidental. Under the
assumptions \fr{cond1} and \fr{cond2}, the optimal (in the minimax
sense) convergence rates in the discrete model are determined by the
total number of observations, $n$, and coincide with the optimal
convergence rates in the continuous model.

In this paper, we consider functional deconvolution in a periodic
setting, that is, we assume that, for fixed $u \in U$, $f(\cdot)$ and
$g(u, \cdot)$ are periodic functions with period on the unit
interval $T$. Note that the periodicity assumption appears naturally
in the above mentioned special models which \fr{convcont} and
\fr{convdis} generalize, and allows one to explore ideas considered
in the above cited papers to the proposed functional deconvolution
framework. Moreover, not only for theoretical reasons but also for
practical convenience [see Johnstone, Kerkyacharian, Picard and
Raimondo (\citeyear{Joh2004}), Sections 2.3, 3.1--3.2], we use band-limited
wavelet bases, and in particular the periodized Meyer wavelet basis
for which fast algorithms exist [see Kolaczyk (\citeyear{Kol1994}) and Donoho and
Raimondo (\citeyear{DonRai2004})].

In what follows, we derive minimax lower bounds for the $L^2$-risk
in models \fr{convcont} and \fr{convdis} when $f(\cdot)$ is assumed
to belong to a Besov ball and $g(\cdot,\cdot)$ is assumed to possess
some smoothness properties, including both regular-smooth and
super-smooth convolutions. Furthermore, we propose an adaptive
wavelet estimator of $f(\cdot)$ and show that this estimator is
asymptotically optimal (in the minimax sense), or near-optimal
within a logarithmic factor, in a wide range of Besov balls. We also
compare models \fr{convcont} and \fr{convdis}, and investigate when
the availability of continuous data gives advantages over
observations at the asymptotically large number of points.

The paper is organized as follows. In Section~\ref{estconstr}, we
describe the construction of a wavelet estimator of $f(\cdot)$
for both the continuous model \fr{convcont} and the discrete model
\fr{convdis}. In Section~\ref{lowbounds}, we derive minimax lower
bounds for the $L^2$-risk, based on observations from either the
continuous model \fr{convcont} or the discrete model \fr{convdis},
when $f(\cdot)$ is assumed to belong to a Besov ball and
$g(\cdot,\cdot)$ is assumed to possess some smoothness properties,
including both regular-smooth and super-smooth convolutions. In
Section~\ref{upbounds}, we demonstrate that the wavelet estimator
derived in Section \ref{estconstr} is adaptive and asymptotically
optimal (in the minimax sense), or near-optimal within a logarithmic
factor, in a wide range of Besov balls. In Section~\ref{applic}, we
discuss particular examples for both continuous and discrete
settings. We conclude in Section~\ref{discrcont} with a discussion
on the interplay between continuous and discrete models. Finally, in
Section~\ref{append}, we provide some auxiliary
statements as well as the proofs of the theoretical results obtained
in the earlier sections.

\section{Construction of a wavelet estimator}
\label{estconstr}

Let $\ph^*(\cdot)$ and $\psi^*(\cdot)$ be the Meyer scaling and
mother wavelet functions, respectively [see, e.g., Meyer (\citeyear{Mey1992}) or
Mallat (\citeyear{Mal1999})]. As usual,
\[
\ph^*_{jk}(x) = 2^{j/2}\ph^*(2^jx-k), \qquad\psi^*_{jk}(x) =
2^{j/2}\psi^*(2^jx-k), \  j,k \in\ints,
\]
are, respectively, the dilated and translated Meyer scaling and
wavelet functions at resolution level $j$ and scale position
$k/2^j$. (Here, and in what follows, $\ints$ refers to the set of
integers.) Similarly to Section 2.3 in Johnstone, Kerkyacharian,
Picard and Raimondo (\citeyear{Joh2004}), we obtain a periodized version of Meyer
wavelet basis by periodizing the basis functions
$\{\ph^*(\cdot),\psi^*(\cdot)\}$, that is,
\[
\ph_{jk}(x) = \sum_{i \in\ints} 2^{j/2} \ph^*\bigl(2^j (x +i) - k\bigr),
\qquad\psi_{jk}(x) = \sum_{i \in\ints} 2^{j/2} \psi^*\bigl(2^j (x +i) -
k\bigr).
\]

In what follows, $\langle\cdot,\cdot\rangle$ denotes the inner
product in the Hilbert space $L^2(T)$ (the space of
squared-integrable functions defined on the unit interval $T$),
that is, $\langle f,g \rangle= \int_{T}f(t)\overline{g(t)} \,dt$ for
$f,g \in L^2(T)$. [Here, and in what follows, $\overline{h(\cdot)}$\vspace*{2pt}
(resp. $\bar{h}$) denotes the conjugate of the complex function
$h(\cdot)$ (resp. complex number $h$); $h(\cdot)$ (resp. $h$) is
real if and only if $\overline{h(\cdot)}=h(\cdot)$ (resp.
$\bar{h}=h$).]

Let $e_m(t) = e^{i 2 \pi m t}$, $m \in\ints$, and, for any (primary
resolution level) $j_0 \geq0$ and any $j \geq j_0$, let
\[
\phmlk= \langle e_m, \ph_{j_0k} \rangle,  \qquad \psimjk= \langle
e_m, \psijk\rangle, \qquad \fm= \langle e_m, f \rangle
\]
be the Fourier coefficients of $\ph_{jk}(\cdot)$, $\psi_{jk}(\cdot)$
and $f(\cdot)$, respectively. Denote
%
\begin{equation}\label{funh}
h(u,t) = \int_T f(x)g(u,t-x)\, dx, \qquad u \in U,  \  t \in T.
\end{equation}
For each $u \in U$, denote the functional Fourier
coefficients by
\begin{eqnarray*}\label{Fouriercont}
\hmu& = &\langle e_m,h (u,\cdot)
\rangle, \qquad    \ymu=
\langle e_m, y(u,\cdot) \rangle, \\
\gmu& = &\langle e_m, g(u,\cdot) \rangle,  \qquad  \zmu= \langle
e_m, z(u,\cdot) \rangle. \nonumber
\end{eqnarray*}

If we have the continuous model \fr{convcont}, then, by using
properties of the Fourier transform, for each $u \in U$, we have
$\hmu= \gmu\fm$ and
%
\begin{equation}\label{finaleq}
\ymu= \gmu\fm+ \frac{1}{\sqrt{n}} \zmu,
\end{equation}
where $\zmu$ are generalized
one-dimensional Gaussian processes such that
%
\begin{equation} \label{zmu}
\EE[z_{m_1}(u_1) z_{m_2}(u_2)] = \delta_{m_1,m_2} \delta(u_1-u_2),
\end{equation}
where $\delta_{m_1,m_2}$ is Kronecker's delta. In order to find the
functional Fourier coefficients $\fm$ of $f(\cdot)$, we multiply
both sides of \fr{finaleq} by $\bgmu$ and integrate over $u \in U$.
The latter yields the following estimators of $\fm$:
%
\begin{equation}\label{fmexprc}
\hfm=\biggl( \iab\bgmu\ymu \,du  \biggr)  \Big/  \biggl( \iab|\gmu|^2 \,du  \biggr).
\end{equation}
[Here, we adopt the convention that when $a=b$
the estimator $\hfm$ takes the form $\hfm= \bgma\yma/|\gma|^2$.]

If we have the discrete model \fr{convdis}, then, by using
properties of the discrete Fourier transform, for each
$l=1,2,\ldots,M$, \fr{finaleq} takes the form
%
\begin{equation}\label{finaleqdis}
\ym(u_l) = \gm (u_l) \fm+ \frac{1}{\sqrt{N}} z_{ml},
\end{equation}
where
$z_{ml}$ are standard Gaussian random variables, independent for
different $m$ and $l$. Similarly to the continuous case, we multiply
both sides of \fr{finaleqdis} by $\overline{g_m(u_l)}$ and add them
together to obtain the following estimators of $\fm$:
%
\begin{equation}\label{fmexprd}
\hfm= \Biggl( \sumM\overline{\gm(u_l)} \ym(u_l)  \Biggr) \bigg/  \Biggl(
\sumM|\gm(u_l)|^2  \Biggr).
\end{equation}
[Here, and in what
follows, we abuse notation and $\fm$ refers to both functional
Fourier coefficients and their discrete counterparts. Note also that
$\ym(u_l)$, $\gm(u_l)$ and $z_{ml}$ are, respectively, the
discrete versions of the functional Fourier coefficients $\ymu$,
$\gmu$ and $\zmu$.]

Note that, using the periodized Meyer wavelet basis described above
and for any $j_0 \geq0$, any (periodic) $f(\cdot) \in L^2(T)$ can
be expanded as
%
\begin{equation}
\label{funf}
f(t) = \sum_{k=0}^{2^{j_0}-1} a_{j_0k} \ph_{j_0k}
(t) + \sum_{j=j_0}^\infty\sum_{k=0}^{2^j -1} b_{jk} \psijk(t).
\end{equation}
Furthermore, by Plancherel's formula, the scaling
coefficients, $\ajk=\langle f, \ph_{j_0k} \rangle$, and the wavelet
coefficients, $\bjk=\langle f,\psi_{jk}\rangle$, of $f(\cdot)$ can
be represented as
%
\begin{equation} \label{alkandblk}
\ajk= \sum_{m \in C_{j_0}} \fm
\overline{\phmlk}, \qquad  \bjk= \sum_{m \in C_j} \fm
\overline{\psimjk},
\end{equation}
where $C_{j_0} = \{m\dvtx\phmlk\neq0 \}$ and, for any $j \geq j_0$, $C_j = \{m:
\psimjk\neq0 \}$, both subsets of $2\pi/3 [-2^{j+2}, -2^j] \cup
[2^j, 2^{j+2}]$, due to the fact that Meyer wavelets are band-limited
[see, e.g., Johnstone, Kerkyacharian, Picard and Raimondo
(\citeyear{Joh2004}), Section 3.1]. We naturally estimate $\ajk$ and $\bjk$ by
substituting $\fm$ in \fr{alkandblk} with \fr{fmexprc} or
\fr{fmexprd}, that is,
%
\begin{equation}\label{coefest}
\hajk= \sum_{m \in C_{j_0}} \hfm
\overline{\phmlk}, \qquad  \hbjk= \sum_{m \in C_j} \hfm
\overline{\psimjk}.
\end{equation}

We now construct a block thresholding wavelet estimator of
$f(\cdot)$. For this purpose, we divide the wavelet coefficients at
each resolution level into blocks of length $\ln n$. Let $\Aj$ and
$\Ujr$ be the following sets of indices:
\begin{eqnarray*}
\Aj &=& \{r \mid r=1,2,\ldots, 2^j/\ln n \},
\\
\Ujr &=& \{
k \mid k = 0,1, \ldots, 2^j-1; (r-1) \ln n \leq k \leq r \ln n -1\}.
\end{eqnarray*}
Denote
%
\begin{equation}
\label{bjr}
\Bjr= \sumku\bjk^2, \qquad  \hBjr= \sumku\hbjk^2.
\end{equation}
Finally, for any $j_0 \geq0$, we reconstruct
$f(\cdot)$ as
%
\begin{equation}\label{fest}
\quad \hat{f}_n(t) = \sum_{k=0}^{2^{j_0} -1} \hajk
\ph_{j_0k} (t) + \sum_{j=j_0}^{J-1} \sumr\sumku\hbjk\II(|\hBjr|
\geq\lambda_{j}) \psijk(t),
\end{equation}
where $\II(A)$ is the
indicator function of the set $A$, and the resolution levels $j_0$
and $J$ and the thresholds $\lambda_{j}$ will be defined in
Section~\ref{upbounds}.

In what follows, we use the symbol $C$ for a generic positive
constant, independent of $n$,  which may take different
values at different places.

\section{Minimax lower bounds for the $L^2$-risk over Besov balls}
\label{lowbounds}

Among the various characterizations of Besov spaces for periodic
functions defined on $L^p(T)$ in terms of wavelet bases, we recall
that for an $r$-regular multiresolution analysis with $0< s < r$ and
for a Besov ball $B_{p,q}^s (A)$ of radius $A>0$ with $1 \leq p,q
\leq\infty$, one has that, with $s' = s+1/2-1/p$,
\begin{eqnarray}
\label{bpqs}
B_{p,q}^s (A) & = & \Biggl\{ f(\cdot) \in L^p(T)\dvtx  \Biggl(
\sum_{k=0}^{2^{j_0}-1}|a_{j_{0}k}|^p
\Biggr)^{1/p}\nonumber\\[-8pt]\\[-8pt]
&&\phantom{\Biggl\{ f(\cdot) \in L^2(T)\dvtx}
{}+  \Biggl( \sum_{j=j_0}^{\infty} 2^{js'q}  \Biggl( \sumk
|\bjk|^p  \Biggr)^{q/p}  \Biggr)^{{1/q}} \leq A  \Biggr\},\nonumber
\end{eqnarray}
with respective sum(s) replaced by maximum if $p=\infty$ or
$q=\infty$ [see, e.g., Johnstone, Kerkyacharian, Picard and Raimondo
(\citeyear{Joh2004}), Section 2.4]. (Note that, for the Meyer wavelet basis,
considered in Section \ref{estconstr}, $r=\infty$.)

We construct below minimax lower bounds for the $L^2$-risk, for both
the continuous model \fr{convcont} and the discrete model
\fr{convdis}. For this purpose, we define the minimax $L^2$-risk
over the set $\Omega$ as
\[
R_n (\Omega) = \inf_{\tilde{f}_n} \, \sup_{f \in\Omega} \EE\|
\tilde{f}_n - f \|^2,
\]
where $\| g \|$ is the $L^2$-norm of a function $g(\cdot)$ and the
infimum is taken over all possible estimators $\tilde{f}_n(\cdot)$ (measurable functions taking their values in a set containing~$\Omega$)
of $f(\cdot)$, based on observations from either the continuous model
\fr{convcont} or the discrete model \fr{convdis}. [Here, and in what
follows, the expectation is taken under the true $f(\cdot)$, and it
is assumed that the function class $\Omega$ contains
$f(\cdot)$.]

In what follows, we shall evaluate a lower bound for $R_n (B_{p,q}^s
(A))$. Denote
\[
s^* = s+1/2-1/p', \qquad p'=\min(p,2),
\]
and, for $\kappa=1,2$, define
%
\begin{equation}\label{taum}
\tau_\kappa(m) = \cases{
\displaystyle\iab|\gmu|^{2\kappa} \,du, &\quad in the continuous case, \cr
\displaystyle\frac{1}{M} \sumM|\gm(u_l)|^{2\kappa}, &\quad in the discrete case.
}
\end{equation}
[Here, we adopt the convention that when $a=b$, $\tau_\kappa
(m)$ takes the form $\tau_\kappa(m) = |g_m(a)|^{2\kappa}$,
$\kappa=1,2$.] Assume that for some constants $\nu\in\RR$, $\af
\geq0$, $\beta>0$ and $K_1 >0$, independent of $m$, the choice of
$M$ and the selection points $u_l$, $l=1,2,\ldots,M$,
%
\begin{equation}\label{cond1}
\tau_1(m)\leq K_1 |m|^{-2\nu} \exp (-\af|m|^\beta ), \qquad\nu>0
\mbox{ if } \af=0.
\end{equation}
[Following Fan
(\citeyear{Fan1991}), we say that the function $g(\cdot,\cdot)$ is
\textit{regular-smooth} if $\af=0$ and is \textit{super-smooth} if $\af
>0$.]

The following statement provides the minimax lower bounds for the
$L^2$-risk.

\begin{thm} \label{th:lower}
Let $\{\phi_{j_0,k}(\cdot),\psi_{j,k}(\cdot)\}$ be the periodic
Meyer wavelet basis discussed in Section \ref{estconstr}. Let $s
>\max(0,1/p-1/2)$, $1 \leq p \leq\infty$, $1 \leq q \leq\infty$ and
$A>0$. Then, under the assumption \fr{cond1}, as $n \rightarrow
\infty$,
%
\begin{equation}\label{low1}
\qquad
R_n (B_{p,q}^s (A)) \geq\cases{
\displaystyle Cn^{-2s/(2s+2\nu+1)}, &\quad if $\af=0$,  $\nu(2-p) < p{s^*}$,
\cr
\displaystyle C  \biggl( \frac{\ln n}{n}  \biggr)^{2{s^*}/(2s^*+2\nu)}, &\quad
if   $  \af=0$,  $\nu(2-p) \geq p{s^*}$,
\cr
\displaystyle C (\ln n)^{-2{s^*}/\beta}, &\quad if  $  \af>0$.
}
\end{equation}
\end{thm}

\begin{remark}
The two different lower bounds for $\af=0$ in (\ref{low1})
refer to the dense case $[\nu(2-p) < p{s^*}]$ when the worst
functions $f(\cdot)$ (i.e., the hardest functions to estimate) are
spread uniformly over the unit interval $T$, and the sparse case
$[\nu(2-p)\geq p{s^*}]$ when the worst functions $f(\cdot)$ have
only one nonvanishing wavelet coefficient. Note also that the
restriction $s>\max(0,1/p-1/2)$, $1\leq p\leq \infty$ and $1\leq q \leq
\infty$ ensures that the corresponding Besov spaces are embedded in
$L^2(T)$.
\end{remark}

\section{Minimax upper bounds for the $L^2$-risk over Besov balls}
\label{upbounds} \setcounter{equation}{0}

Recall (\ref{cond1}) from Section \ref{lowbounds}, and assume
further that for the constants $\nu\in\RR$, $\af\geq0$ and
$\beta>0$, and for a constant $K_2 >0$, independent of $m$, the
choice of $M$ and the selection points $u_l$, $l=1,2,\ldots,M$, with
$K_2 \leq K_1$,
%
\begin{equation}
\label{cond2}
\tau_1 (m) \geq K_2 |m|^{-2\nu} \exp (-\af
|m|^\beta ), \qquad\nu>0 \mbox{ if } \af=0.
\end{equation}
For any $j \geq j_0$, let $|C_j|$ be the
cardinality of the set $C_j$; note that, for Meyer wavelets, $|C_j|
= 4 \pi2^j $ [see, e.g., Johnstone, Kerkyacharian, Picard and
Raimondo (\citeyear{Joh2004}), page 565]. Let also
%
\begin{equation}\label{Delt}
\Delta_\kappa(j) =
\frac{1}{|C_j|} \sum_{m \in C_j} \tau_{\kappa} (m) [\tau_1 (m)]^{-2
\kappa}, \qquad \kappa=1,2.
\end{equation}
Then, direct calculations
yield that
%
\begin{equation}\label{delta1}
\Delta_1 (j) \leq\cases{
\displaystyle c_1 2^{2 \nu j}, &\quad
if  $  \af=0$, \cr 
\displaystyle c_2 2^{2 \nu j} \exp \biggl( \af \biggl( \frac{8\pi}{3}  \biggr) ^\beta
2^{j \beta}  \biggr), &\quad if   $  \af> 0$. 
}
\end{equation}
[Note that since the functional Fourier coefficients
$g_m(\cdot)$ are known, the positive constants $c_1$ and $c_2$ in
\fr{delta1} can be evaluated explicitly.]

Consider now the two cases $\af=0$ (regular-smooth) and $\af>0$
(super-smooth) separately. Choose $\jo$ and $J$ such that
%
\begin{eqnarray}\label{jpower}
2^\jo &=& \ln n, \qquad\hspace*{45pt}  2^J = n^{1/(2\nu+1)} \qquad \mbox{if}\    \af
=0. \\
2^\jo &=& \frac{3}{8 \pi} \biggl(\frac{\ln n}{2\af}
\biggr)^{1/\beta},\qquad   2^J=2^{\jo}, \qquad\hspace*{21pt} \mbox{if}\    \af>0.
\label{jexp}
\end{eqnarray}
[Since $\jo>J-1$ when $\af>0$, the wavelet
estimator \fr{fest} only consists of the first (linear) part and,
hence, $\lambda_j$ does not need to be selected in this case.] Set, for
some positive constant $d$,
%
\begin{equation} \label{lamj}
\lambda_j = d  n^{-1} 2^{2 \nu j} \ln
n\qquad     \mbox{if }  \af=0.
\end{equation}
Note that the choices
of $\jo$, $J$ and $\lambda_j$ are independent of the parameters, $s$,
$p$, $q$ and $A$ (that are usually unknown in practical situations)
of the Besov ball $B_{p,q}^s (A)$; hence, the wavelet estimator
\fr{fest} is adaptive with respect to these parameters.

The proof of the minimax upper bounds for the $L^2$-risk is based on
the following two lemmas.

\begin{lemma} \label{l:coef}
Let the assumption \fr{cond2} be valid, and let the estimators
$\hajk$ and $\hbjk$ of the scaling and wavelet coefficients $\ajk$
and $\bjk$, respectively, be given by the formula \fr{coefest} with
$\hfm$ defined by \fr{fmexprc} in the continuous model and by
\fr{fmexprd} in the discrete model. Then, for $\kappa=1,2$, and for
all $j \geq j_0$,
\begin{eqnarray}
\EE|\hajk- \ajk|^2 & \leq & C n^{-1} \Delta_1 (\jo), \label{ha}\\
\EE|\hbjk- \bjk|^{2 \kappa} & \leq & C n^{-\kappa} \Delta_\kappa
(j). \label{hb}
\end{eqnarray}
Moreover, under the assumptions
\fr{cond1} and \fr{cond2} with $\alpha=0$, for all $j \geq j_0$,
%
\begin{equation} \Delta_2 (j) \leq C 2^{4(2\nu-\nu_1)j},
\label{fanD2}
\end{equation}
for any $0 < \nu_1 \leq\nu$.
\end{lemma}

\begin{lemma} \label{l:deviation}
Let the estimators $\hbjk$ of the wavelet coefficients $\bjk$ be
given by the formula \fr{coefest} with $\hfm$ defined by
\fr{fmexprc} in the continuous model and by \fr{fmexprd} in the
discrete model. If $\mu$ is a positive constant large enough and
$\af=0$ in the assumption \fr{cond2}, then, for all $j \geq j_0$,
%
\begin{equation}\label{probbound}
\qquad \PP \Biggl( \sumku|\hbjk- \bjk|^2 \geq0.25 \mu^2 n^{-1}
2^{2\nu j} \ln n  \Biggr) \leq n^{-(8\nu- 4 \nu_1
+2)/(2\nu+1)},
\end{equation}
for any $0 < \nu_1 \leq\nu$.
\end{lemma}

Lemmas \ref{l:coef} and \ref{l:deviation} allow to state the
following minimax upper bounds for the $L^2$-risk of the wavelet
estimator $\hat{f}_n(\cdot)$ defined by \fr{fest}, with $\jo$ and
$J$ given by \fr{jpower} (if $\af= 0$) or \fr{jexp} (if $\af
>0$). Set $(x)_{+} = \max(0,x)$, and define
%
\begin{equation}
\ro_1 = \cases{
\dfrac{(2\nu+1)(2-p)_+}{p(2s+2\nu+1)}, &\quad if $\nu(2-p) <
p{s^*}$, \cr 
\dfrac{(q-p)_+}{q}, & \quad if  $   \nu(2-p) = p{s^*}$,
\label{rovalue}\cr
0, & \quad if  $  \nu(2-p) > p{s^*}$. 
}
\end{equation}

\begin{thm} \label{th:upper}
Let $\hat{f}_n(\cdot)$ be the wavelet estimator defined by
\fr{fest}, with $\jo$ and $J$ given by \fr{jpower} (if $\af= 0$) or
\fr{jexp} (if $\af>0$). Let $s > 1/p'$, $1 \leq p \leq\infty$, $1
\leq q \leq\infty$ and $A>0$. Then, under the assumption
\fr{cond2}, as $n \rightarrow\infty$,
%
\begin{equation}\label{up}
\sup_{f \in B_{p,q}^s
(A)} \EE\|\hf-f\|^2 \leq\cases{
Cn^{-2s/(2s+2\nu+1)} (\ln n)^{\ro_1}, \cr
\qquad
\mbox{if }\af=0, \nu(2-p) < p{s^*},
\cr 
C  \biggl( \dfrac{\ln n}{n}  \biggr)^{2{s^*}/(2s^*+2\nu)}  (
\ln n )^{\ro_1}, \cr
\qquad \mbox{if  } \af=0, \nu(2-p) \geq p{s^*}, \cr
C (\ln n)^{-2{s^*}/\beta}, \cr
\qquad\mbox{if }\af>0.
}
\end{equation}
\end{thm}

\begin{remark}
In the discrete model, assumptions \fr{cond1} and \fr{cond2}
require the value of $\tau_1(m)$ to be independent of the choice of
$M$ and the selection of points $u_l$, $l=1,2, \ldots, M$. If
assumptions \fr{cond1} and \fr{cond2} hold, then the minimax
convergence rates in discrete and continuous models coincide and are
independent of the configuration of the points $u_l$, $l=1,2,
\ldots, M$. Moreover, the wavelet estimator \fr{fest} is
asymptotically optimal (in the minimax sense) no matter what the
value of $M$ is. It is quite possible, however, that in the discrete
model, conditions \fr{cond1} and \fr{cond2} both hold but with
different values of $\nu$, $\af$ and $\beta$. In this case, the
upper bounds for the risk in the discrete model may not coincide with
the lower bounds and with the minimax convergence rates in the
continuous model. Proposition \ref{th:condisrates} in Section \ref{discrcont} provides sufficient
conditions for the minimax convergence rates in discrete and
continuous models to coincide and to be independent of $M$ and the
configuration of the points $u_l$, $l=1,2, \ldots, M$. These
conditions also guarantee asymptotical optimality of the wavelet
estimator \fr{fest}, and can be viewed as some kind of uniformity
conditions. If conditions of Proposition \ref{th:condisrates} are violated, then the
rates of convergence in the discrete model depend on the choice of
$M$ and $u_l$, $l=1,2, \ldots, M$, and some recommendations on their
selection should be given. Furthermore, optimality issues become
much more complex when $\tau_1(m)$ is not uniformly bounded from
above and below (see the discussion in Section \ref{discrcont}).
\end{remark}

\begin{remark}
Theorems \ref{th:lower} and \ref{th:upper} imply that, for the
$L^2$-risk, the wavelet estimator $\hat{f}_n(\cdot)$ defined by
\fr{fest} is asymptotically optimal (in the minimax sense), or
near-optimal within a logarithmic factor, over a wide range of Besov
balls $B_{p,q}^s (A)$ of radius $A>0$ with $s
> \max(1/p, 1/2)$, $1 \leq p \leq\infty$ and $1 \leq q \leq\infty
$. In
particular, in the cases when (1) $\af> 0$, (2) $\af=0$, $\nu(2-p) <
p{s^*}$ and $2 \leq p \leq\infty$, (3)~$\af=0$, $\nu(2-p)
> p{s^*}$, and (4) $\af=0$, $\nu(2-p) = p{s^*}$ and $1 \leq q \leq p$, the
estimator \fr{fest} is asymptotically optimal (lower and upper
bounds coincide up to a multiplicative constant), that is,
\[
R_n (B_{p,q}^s (A)) \asymp\cases{
n^{-2s/(2s+2\nu+1)}, & \quad if  $\af=0$,  $\nu(2-p) <
p{s^*}$, $2\leq p \leq\infty$,\cr
\biggl( \dfrac{\ln n}{n}  \biggr)^{2{s^*}/(2s^*+2\nu)}, & \quad
if  $\af=0$,  $\nu(2-p) > p{s^*}$, \cr
&\quad or    $\af=0$,  $\nu(2-p) = p{s^*}$, $1 \leq q \leq p$, \cr
(\ln n)^{-2{s^*}/\beta}, & \quad if $   \af>0$.
}
\]
On the other hand, in the case when $\af=0$, $\nu(2-p) < p{s^*}$
and $1 \leq p<2$ or $\af=0$, $\nu(2-p) = p{s^*}$ and $1 \leq p<q$,
the wavelet estimator $\hat{f}_n(\cdot)$ defined by \fr{fest} is
asymptotically near-optimal within a logarithmic factor, that is,
\[\label{low}
\sup_{f \in B_{p,q}^s (A)} \EE\|\hf-f\|^2 \leq\cases{
C n^{-2s/(2s+2\nu+1)} (\ln
n)^{(2\nu+1)(2-p)/(p(2s+2\nu+1))}, \cr\qquad    \mbox{if }
\af=0,\    \nu(2-p) <
p{s^*}, \cr
\qquad 1 \leq p < 2,\cr
C  \biggl( \dfrac{\ln n}{n}  \biggr)^{2{s^*}/(2s^*+2\nu)} (\ln
n)^{(1-p/q)},  \cr
\qquad\mbox{if }\af=0, \nu(2-p) =
p{s^*},\cr
\qquad 1 \leq p < q.
}
\]
[Here, and in what follows, $g_1(n) \asymp g_2(n)$ denotes $0 <
\liminf(g_1(n)/g_2(n)) \leq\limsup(g_1(n)/g_2(n)) < \infty$ as $n
\rightarrow\infty$.]
\end{remark}

\begin{remark}
For the $L^2$-risk, the upper bounds \fr{up} are tighter than
those obtained by Chesneau (\citeyear{Che2008}) for the regular-smooth case [i.e.,
$\alpha=0$ in \fr{cond1} and \fr{cond2}] in the case of the
standard deconvolution model [i.e., when $a=b$ in (\ref{convcont})],
although the difference is only in the logarithmic factors. More
specifically, the following minimax upper bounds obtained in
Chesneau (\citeyear{Che2008}) for the $L^2$-risk, as $n \rightarrow\infty$:
\begin{equation}\label{eq:chesFFF}
\sup_{f \in B_{p,q}^s (A)} \EE\|\hf-f\|^2 \leq\cases{
Cn^{-2s/(2s+2\nu+1)} (\ln n)^{\ro_2}, \cr
\qquad \mbox{if }\af=0, \nu(2-p) < p{s'},
\cr
C  \biggl( \dfrac{\ln n}{n}  \biggr)^{2{s'}/(2s'+2\nu)}  (\ln
n )^{\ro_2}, \cr \qquad \mbox{if }   \af=0, \nu(2-p) \geq p{s'},
}
\end{equation}
where
%
\begin{equation}\label{rostarvalue}
\ro_2 = \cases{
\dfrac{2s \II(1 \leq p<2)}{ 2s+2\nu+1}, &\quad if $ \nu(2-p)
< p{s'}$, \cr
\dfrac{(2q-p)_+}{q}, &\quad if $   \nu(2-p) = p{s'}$, \cr
0, &\quad if $\nu(2-p) > p{s'}$.
}
\end{equation}
[Here, and in what follows, $\II(A)$ is the indicator function
of the set $A$.] Note that when $2\leq p \leq\infty$, $s^*=s \leq
s'$, and only the dense case appears; hence, in this case, the dense
cases and the corresponding convergence rates in the minimax upper
bounds given by \fr{rovalue}--\fr{up} and
\fr{eq:chesFFF}--\fr{rostarvalue} coincide since $\nu(2-p) <
p{s^*}=ps<ps'$. On the other hand, when $1 \leq p <2$, $s^*=s'$,
both the dense and sparse cases appear; hence, in this case, both
the dense and sparse cases and the corresponding convergence rates
in the minimax upper bounds given by \fr{rovalue}--\fr{up} and
\fr{eq:chesFFF}--\fr{rostarvalue} coincide. Looking now at
\fr{rovalue} and \fr{rostarvalue}, we see that $\ro_2=\ro_1$ only
when $\nu(2-p)
> p{s'}$. On the other hand, $\ro_2 > \ro_1$ when $1 \leq p<2$ and $
\nu(2-p) < p{s'}$ since $(2/p -1)(2\nu+1) - 2s = 2(2\nu- p \nu- p
s')/p <0$, and it is obvious that $\ro_2 > \ro_1$ when $\nu(2-p) =
p{s'}$. However, we believe that the slight superiority in the
minimax convergence rates for the $L^2$-risk obtained in Theorems
\ref{th:lower} and \ref{th:upper} is due not to a different
construction of the wavelet estimator but to a somewhat different
way of evaluating the minimax upper bounds.
\end{remark}

\begin{remark}
Unlike Chesneau (\citeyear{Che2008}) who only considered minimax upper bounds
for the regular-smooth case [i.e., $\alpha=0$ in \fr{cond1} and
\fr{cond2}] in the standard deconvolution model [i.e., when $a=b$ in
(\ref{convcont})], Theorems \ref{th:lower} and \ref{th:upper}
provide minimax lower and upper bounds (in the $L^2$-risk) for both
regular-smooth and super-smooth convolutions [i.e., $\alpha>0$ in
\fr{cond1} and \fr{cond2}], not only for the standard deconvolution
model but also for its discrete counterpart [i.e., when $M=1$ in
(\ref{convdis})].
\end{remark}

\begin{remark}
The wavelet estimator $\hat{f}_n(\cdot)$ defined by \fr{fest}
is adaptive with respect to the unknown parameters $s$, $p$, $q$ and
$A$ of the Besov ball $B_{p,q}^s (A)$ but is not adaptive with
respect to the parameters $\alpha$, $\beta$ and $\nu$ in \fr{cond1}
and \fr{cond2}. It seems that it is impossible to achieve adaptivity
with respect to $\beta$ in the super-smooth case ($\alpha
>0$) because of the very fast exponential growth of the variance.
However, in the regular-smooth case ($\af=0$), one can construct a
wavelet estimator which is adaptive with respect to the unknown
parameter $\nu$. Choose $\jo$ and $J$ such that $2^\jo= \ln n$ and
$2^J = n$, and set $ \lambda_j = d^* n^{-1} \ln n  \Delta_1 (j)$,
where $d^*$ is large enough. Note that $\Delta_1 (j)$ can be
calculated whenever the functional Fourier coefficients $\gm(\cdot)$
are available. Also, $K_1^* 2^{2 \nu j} \leq\Delta_1 (j) \leq K_2^*
2^{2 \nu j}$ for some positive constants $K_1^*$ and $K_2^*$ which
depend on the particular values of the constants in the conditions
\fr{cond1} and \fr{cond2}. Therefore, in this situation, by
repeating the proof of Theorem \ref{th:upper} with these new values of the
parameters involved, one can easily verify that the optimal
convergence rates in Theorem \ref{th:upper} still hold as long as $d^*$ is large
enough. How large should be ``large enough''? Direct calculations
show that $d^*$ should be such that $(0.5 d^*
\sqrt{K_2^*/K_1^*}-1)^2 \geq8\nu+2$. Since $K_1^*$, $K_2^*$ and
$\nu$ are unknown, it is impossible to evaluate the lower bound for
$d^*$. However, one can replace $d^*$ by a slow-growing function of
$n$, say $\ln\ln n$, leading to, at most, an extra $\ln\ln n$
factor in the obtained maximal $L^2$-risk.
\end{remark}

\begin{remark}
We finally note that, although we have only considered
$L^2$-risks in our analysis, the results obtained in Theorems
\ref{th:lower} and \ref{th:upper} can be extended to the case of
$L^{\pi}$-risks ($1 \leq\pi< \infty$). Analogous statements  to the
ones given in Theorems \ref{th:lower} and \ref{th:upper} but for a
wider variety of risk functions can be obtained using the
unconditionality and Temlyakov properties of Meyer wavelets [see,
e.g., Johnstone, Kerkyacharian, Picard and Raimondo (\citeyear{Joh2004}),
Appendices A and B]. The details in the derivation of these statements
should, however, be carefully addressed.
\end{remark}

\section{Examples in continuous and discrete models}
\label{applic}

The functional deconvolution model (\ref{convcont}) can be viewed as
a generalization of a multitude of inverse problems in mathematical
physics where one needs to recover initial or boundary conditions on
the basis of observations of a noisy solution of a partial
differential equation. Lattes and Lions (\citeyear{LatLio1967}) initiated research in
the problem of recovering the initial condition for parabolic
equations based on observations in a fixed-time strip. This problem
and the problem of recovering the boundary condition for elliptic
equations based on observations in an internal domain were studied
in Golubev and Khasminskii (\citeyear{GolKha1999}). More specifically, by studying
separately the heat conductivity equation or the Laplace equation on
the unit circle, and assuming that the unknown initial or boundary
condition belongs to a Sobolev ball, Golubev and Khasminskii (\citeyear{GolKha1999})
obtained some \textit{linear} and \textit{nonadaptive} solutions to
the particular problem at hand; see also Golubev (\citeyear{Gol2004}) for a linear
adaptive estimator for the Laplace equation on the circle based on
the principle of minimization of penalized empirical risk. We also
note that, unlike Golubev and Khasminskii (\citeyear{GolKha1999}) and Golubev (\citeyear{Gol2004})
who considered sharp asymptotics, we focus our study on rate
optimality results. [Note that the estimation of the unknown initial
condition for the heat conductivity equation, allowing also for
missing data, has been recently considered by Hesse (\citeyear{Hes2007});
however, this latter paper deals with the density deconvolution model and the
approach given therein varies from the approach of Golubev and
Khasminskii (\citeyear{GolKha1999}) and Golubev (\citeyear{Gol2004}), and it seems to be having a
different agenda.]

In view of the general framework developed in this paper, however,
the inverse problems mentioned above can all be expressed as a
functional deconvolution problem, so that all techniques studied in
Sections~\ref{estconstr}--\ref{upbounds} can be directly applied, to
obtain \textit{linear}/\textit{nonlinear} and \textit{adaptive} solutions
over a wide range of Besov balls. Such solutions are provided in
Examples \ref{ex1}--\ref{exam4} below which discuss some of the most common inverse
problems in mathematical physics which have already been studied as
well as some other problems which, to the best of our knowledge,
have not yet been addressed.

On the other hand, in the case when the functional deconvolution
model (\ref{convcont}) is observed at a finite number of distinct
points [see (\ref{convdis})], it can also be viewed as a
multichannel deconvolution model studied in De Canditiis and Pensky
(\citeyear{DeCanPen2004}, \citeyear{DeCanPen2006}). Example \ref{ex5} below deals with this model, providing the
minimax convergence rates (in the $L^2$-risk) for regular-smooth
[i.e., $\alpha=0$ in (\ref{cond1}) and (\ref{cond2})] and
super-smooth [i.e., $\alpha>0$ in (\ref{cond1}) and (\ref{cond2})]
convolutions, and also discussing the case when $M$ can increase
together with $N$; both of these aspects were lacking from the
theoretical analysis described in De Canditiis and Pensky (\citeyear{DeCanPen2006}).

\begin{example}[(\textit{Estimation of the initial condition in the heat
conductivity equation})]\label{ex1}
Let $h(t,x)$ be a solution of the heat conductivity equation
\[
\frac{\partial h (t,x)}{\partial t} = \frac{\partial^2 h
(t,x)}{\partial x^2},\qquad  x \in[0,1],\   t \in[a,b], \   a>0,\
b < \infty,
\]
with initial condition $h(0,x) = f(x)$ and periodic boundary
conditions
\[
h(t,0) = h(t,1), \qquad  \frac{\partial h(t,x)}{\partial x} \bigg|_{x=0}
= \frac{\partial h(t,x)}{\partial x} \bigg|_{x=1}.
\]

We assume that a noisy solution $y(t,x) = h(t,x) + n^{-1/2}z(t,x)$
is observed, where $z(t,x)$ is a generalized two-dimensional
Gaussian field with covariance function $\EE[z (t_1, x_1) z (t_2,
x_2)] = \delta(t_1-t_2) \delta(x_1-x_2)$, and the goal is to recover
the initial condition $f(\cdot)$ on the basis of observations
$y(t,x)$. This problem was considered by Lattes and Lions (\citeyear{LatLio1967}) and
Golubev and Khasminskii (\citeyear{GolKha1999}).

It is well known [see, e.g., Strauss (\citeyear{Str1992}), page 48] that, in a
periodic setting, the solution $h(t,x)$ can be written as
%
\begin{equation}\label{heat}
h(t,x) = (4 \pi t)^{-1/2}  \int_0^1 \sum_{k \in\ints}
\exp \biggl\{-\frac{(x +k -z)^2}{4t} \biggr\} f(z)\, dz.
\end{equation}
It is easy to see that \fr{heat} coincides with \fr{funh} with $t$
and $x$ replaced by $u$ and $t$, respectively, and that
\[
g(u,t) = (4 \pi u)^{-1/2} \sum_{k \in\ints} \exp \biggl\{-\frac{(t
+k)^2}{4u} \biggr\}.
\]
Applying the theory developed in
Sections~\ref{estconstr}--\ref{upbounds}, we obtain functional
Fourier coefficients $g_m(\cdot)$ satisfying $\gmu= \exp(- 4 \pi^2
m^2 u)$, and
\[
\tau_1(m)=\iab|\gmu|^2 \,du = C m^{-2} \exp (-8 \pi^2 m^2 a )
\bigl(1 + o(1)\bigr), \qquad  |m| \rightarrow\infty,
\]
so that $\nu=1$, $\af= 8 \pi^2 a$ and $\beta=2$ in both
\fr{cond1} and \fr{cond2}.

Hence, one can construct an adaptive wavelet estimator of the form
\fr{fest}, with $\jo$ and $J$ given by \fr{jexp}, which achieves
minimax (in the $L^2$-risk) convergence rate of order $(\ln
n)^{-\sstar}$ over Besov balls $\Bpqs(A)$ of radius $A>0$ with $s >
\max(1/p, 1/2)$, $1 \leq p \leq\infty$ and $1 \leq q \leq\infty$.
\end{example}

\begin{example}[(\textit{Estimation of the boundary condition for the
Dirichlet problem of the Laplacian on the unit circle})]\label{ex2}
Let $h(x,w)$ be a solution of the Dirichlet problem of the Laplacian
on a region $D$ on the plane
%
\begin{equation} \label{laplace}
\frac{\partial^2 h (x,w)}{\partial
x^2} + \frac{\partial^2 h (x,w)}{\partial w^2} =0,\qquad   (x, w) \in D
\subseteq\RR^2,
\end{equation}
with a boundary $\partial D$
and boundary condition
%
\begin{equation}
\label{eqFFFboun}
h(x,w) |_{\partial D} = F(x,w).
\end{equation}
Consider the situation when $D$ is the unit
circle. Then, it is advantageous to rewrite the function
$h(\cdot,\cdot)$ in polar coordinates as $h(x,w) = h(u,t)$, where $u
\in[0,1]$ is the polar radius and $t \in[0, 2\pi]$ is the polar
angle. Then, the boundary condition in \fr{eqFFFboun} can be
presented as $h(1, t) = f(t)$, and $h(u,\cdot)$ and $f(\cdot)$ are
periodic functions of $t$ with period $2\pi$.

Suppose that only a noisy version $y(u,t) = h(u,t) + n^{-1/2}z(u,t)$
is observed, where $z(u,t)$ is as in Example \ref{ex1}, and that
observations are available only on the interior of the unit circle
with $u \in[0, r_0]$, $r_0 <1$, that is, $a=0$, $b=r_0<1$. The goal is
to recover the boundary condition $f(\cdot)$ on the basis of
observations $y(u,t)$. This problem was investigated in Golubev and
Khasminskii (\citeyear{GolKha1999}) and Golubev (\citeyear{Gol2004}).

It is well known [see, e.g., Strauss (\citeyear{Str1992}), page 161] that the
solution $h(u,t)$ can be written as
\[
h(u,t) = \frac{(1 - u^2)}{2\pi}  \int_0^{2\pi} \frac{f(x)}{1 - 2u
\cos(t-x) + u^2}\, dx.
\]
Applying the theory developed in
Sections~\ref{estconstr}--\ref{upbounds} with $e_m(t) = e^{i m t}$
and
\[
g(u,t) = \frac{1 - u^2}{1 - 2u \cos(t) + u^2},
\]
we obtain functional Fourier coefficients $g_m(\cdot)$ satisfying
$\gmu= C u^{m}$, and
\[
\tau_1(m)=\int_0^{r_0} |\gmu|^2 \,du = C
\exp\{-2 \ln(1/r_0)|m|\},
\]
so that $\nu=0$, $\af= 2 \ln(1/r_0)$ and $\beta=1$ in both
\fr{cond1} and \fr{cond2}.

Hence, one can construct an adaptive wavelet estimator of the form
\fr{fest}, with $\jo$ and $J$ given by \fr{jexp}, which achieves
minimax (in the $L^2$-risk) convergence rate of order $(\ln
n)^{-2\sstar}$ over Besov balls $\Bpqs(A)$ of radius $A>0$ with $s
> \max(1/p, 1/2)$, $1 \leq p \leq\infty$ and $1 \leq q \leq
\infty$.
\end{example}

\begin{example}[(\textit{Estimation of the boundary condition for the
Dirichlet problem of the Laplacian on a rectangle})]\label{ex3}
Consider the problem \fr{laplace}--\fr{eqFFFboun} in Example~\ref{ex2}
above, with the region $D$ being now a rectangle, that is, $(x,w) \in
[0,1] \times[a,b]$, $a>0$, $b< \infty$, and periodic boundary
conditions
\[
h(x,0) = f(x),\qquad  h(0,w) = h(1,w).
\]
Again, suppose that only a noisy version $y(x,w) = h(x,w) +
n^{-1/2}z(x,w)$ is observed, where $z(x,w)$ is as in Example \ref{ex1}, for
$x \in[0,1]$, $w \in[a,b]$, and the goal is to recover the
boundary condition $f(\cdot)$ on the basis of observations $y(x,w)$.

It is well known [see, e.g., Strauss (\citeyear{Str1992}), pages 188, 407] that,
in a periodic setting, the solution $h(x,w)$ can be written as
%
\begin{equation}
h(x,w) = \pi^{-1} \int_0^1\sum_{k \in\ints} \frac{w}{w^2 +
(x+k-z)^2} f(z) \,dz. \label{fanisEx3}
\end{equation}
It is easy to see that \fr{fanisEx3} coincides with \fr{funh} with
$x$ and $w$ replaced by $t$ and $u$, respectively, and that
\[
g(u,t) = \pi^{-1} \sum_{k \in\ints} \frac{u}{u^2+(t+k)^2}.
\]
Applying the theory developed in
Sections~\ref{estconstr}--\ref{upbounds}, we obtain functional
Fourier coefficients $g_m(\cdot)$ satisfying $\gmu= \exp(-2\pi m
u)$, and
\[
\tau_1(m)=\iab|\gmu|^2 \,du = C |m|^{-1} \exp (-4 \pi|m| a )
\bigl(1 + o(1)\bigr),\qquad   |m| \rightarrow\infty,
\]
so that $\nu= 1/2$, $\af= 4 \pi a$ and $\beta=1$ in both
\fr{cond1} and \fr{cond2}.

Hence, one can construct an adaptive wavelet estimator of the form
\fr{fest}, with $\jo$ and $J$ given by \fr{jexp}, which achieves
minimax (in the $L^2$-risk) convergence rate of order $(\ln
n)^{-2\sstar}$ over Besov balls $\Bpqs(A)$ of radius $A>0$ with $s
> \max(1/p, 1/2)$, $1 \leq p \leq\infty$ and $1 \leq q \leq
\infty$.
\end{example}

\begin{example}[(\textit{Estimation of the speed of a wave on a finite
interval})]\label{exam4}
Let $h(t,x)$ be a solution of the wave equation
\[\label{waves}
\frac{\partial^2 h (t,x)}{\partial t^2} = \frac{\partial^2 h
(t,x)}{\partial x^2}
\]
with
initial--boundary conditions
\[
h(0,x)=0, \qquad   \frac{\partial h(t,x)}{\partial t}  \bigg|_{t=0} =
f(x),  \qquad  h(t,0) = h(t,1)=0.
\]
Here, $f(\cdot)$ is a function defined on the unit interval
$[0,1]$, and the objective is to recover $f(\cdot)$ on the basis of
observing a noisy solution $y(t,x) = h(t,x) + n^{-1/2}z(t,x)$, where
$z(t,x)$ is as in Example \ref{ex1}, with $t \in[a,b]$, $a
>0$, $b<1$.

Extending $f(\cdot)$ periodically over the real line, it is
well known that the solution $h(t,x)$ can then be recovered as [see,
e.g., Strauss (\citeyear{Str1992}), page 61]
%
\begin{equation}\label{wav1}
h(t,x) = \tfrac{1}{2} \int_0^1
\II(|x-z| < t) f(z) \,dz,
\end{equation}
so that \fr{wav1} is of the
form \fr{funh} with $g(u,x)= 0.5\II(|x| < u)$ (a boxcar-like
kernel for each fixed $u$), where $u$ in \fr{funh} is replaced by
$t$ in \fr{wav1}. Applying the theory developed in
Sections~\ref{estconstr}--\ref{upbounds}, with $t$ and $x$ replaced
by $u$ and $t$, respectively, we obtain functional Fourier
coefficients $g_m(\cdot)$ satisfying $\gmu= \sin(2 \pi m u)/(2\pi
m)$, and
\begin{eqnarray} \label{ex4}
\tau_1(m) &=& \iab|\gmu|^2\, du\nonumber\\[-8pt]\\[-8pt]
& = & \frac{1}{4\pi^2 m^2}
 \biggl(\frac{b-a}{2} + \frac{\sin(4\pi m a) - \sin(4\pi m b)}{8 \pi m}
\biggr).\nonumber
\end{eqnarray}
Observe that the integral in \fr{ex4} is
always positive, bounded from above by $C m^{-2}$ and from below by
$Cm^{-2} [(b-a) - (2 \pi m)^{-1}]$, so that $\nu=1$ and $\af=0$ in
both \fr{cond1} and \fr{cond2}.

Hence, one can construct an adaptive block thresholding wavelet
estimator of the form \fr{fest}, with $\jo$ and $J$ given by
\fr{jpower}, which achieves the following minimax upper bounds (in
the $L^2$-risk):
\[
\sup_{f \in B_{p,q}^s (A)} \EE\|\hf-f\|^2 \leq\cases{
C n^{-2s/(2s+3)} (\ln n)^{\ro_1}, &\quad if
$s>3(1/p - 1/2)$,\cr
C  \biggl( \dfrac{\ln n}{n}  \biggr)^{s'/(s'+1)}  (\ln n
)^{\ro_1}, &\quad if $s \leq3(1/p - 1/2)$,
}
\]
over Besov balls $\Bpqs(A)$ of radius $A>0$ with $s > 1/p'$, $1
\leq p \leq\infty$ and $1 \leq q \leq\infty$, where $\ro_1 =
3(2/p-1)_+/(2s+3)$ if $s >3(1/p - 1/2)$, $\ro_1 = (1 -p/q)_+$ if $s
= 3(1/p - 1/2)$ and $\ro_1 =0$ if $s < 3(1/p - 1/2)$. [The minimax
lower bounds (in the $L^2$-risk) have the same form with $\ro_1=0$.]
\end{example}

\begin{example}[(\textit{Estimation in the multichannel deconvolution
problem})]\label{ex5}
Consider the problem of recovering $f(\cdot) \in L^2(T)$
on the basis of observing the following noisy convolutions with
known blurring functions $g_l(\cdot)$
%
\begin{equation}\label{multchan}
\qquad Y_{l}(dt) = f*g_l (t) \,dt +
\frac{\sigma_l}{\sqrt{ n}} W_l (dt), \qquad  t \in T,
l=1,2,\ldots,M.
\end{equation}
Here, $\sigma_l$ are known positive
constants and $W_l(t)$ are independent standard Wiener processes.

The problem of considering systems of convolution equations was
first considered by Casey and Walnut (\citeyear{CasWal1994}) in order to evade the
ill-posedness of the standard deconvolution problem, and was adapted
for statistical use (in the density deconvolution model) by Pensky
and Zayed (\citeyear{PenZay2002}). Wavelet solutions to the problem \fr{multchan} were
investigated by De Canditiis and Pensky (\citeyear{DeCanPen2004}, \citeyear{DeCanPen2006}).

Note that deconvolution is the common problem in many areas of
signal and image processing which include, for instance, LIDAR
(Light Detection and Ranging) remote sensing and reconstruction of
blurred images. LIDAR is a laser device which emits pulses,
reflections of which are gathered by a telescope aligned with the
laser [see, e.g., Park, Dho and Kong (\citeyear{JeParWhoDhoJinKon1997}) and Harsdorf
and Reuter (\citeyear{HarReu2000})]. The return signal is used to determine distance
and the position of the reflecting material. However, if the system
response function of the LIDAR is longer than the time resolution
interval, then the measured LIDAR signal is blurred and the
effective accuracy of the LIDAR decreases. If $M$ ($M\geq2$) LIDAR
devices are used to recover a signal, then we talk about a
multichannel deconvolution problem. Note that a discretization of
\fr{multchan} (with $\sigma_l=1$ for $l=1,2,\ldots,M$) leads to the
discrete setup \fr{convdis}.

Adaptive term by term wavelet thresholding estimators for the model
\fr{multchan} were constructed in De Canditiis and Pensky (\citeyear{DeCanPen2006}) for
regular-smooth convolutions [i.e., $\alpha=0$ in (\ref{cond1}) and
(\ref{cond2})]. However, minimax lower and upper bounds were not
obtained by these authors who concentrated instead on upper bounds
(in the $L^{\pi}$-risk, $1 < \pi< \infty$) for the error, for a
fixed response function. Moreover, the case of super-smooth
convolutions [i.e., $\alpha>0$ in (\ref{cond1}) and (\ref{cond2})]
and the case when $M \rightarrow\infty$ have not been treated in De
Canditiis and Pensky (\citeyear{DeCanPen2006}).

Let us now discuss the regular-smooth convolution case treated in De
Canditiis and Pensky (\citeyear{DeCanPen2006}), that is, the case when (in our notation)
$|\gm(u_l)| \sim C_l |m|^{- \nu_l}$ with $0 < C_l <\infty$,
$l=1,2,\ldots,M$. If $M$ is fixed, then
\[
C_* M^{-1} m^{-2\nu_{\min}} \leq\tau_1 (m) \leq C^* m^{-2\nu
_{\min}},
\]
where $\nu_{\min} = \min\{\nu_1, \nu_2,\ldots,\nu_M\}$ and $0<C_*
\leq C_l \leq C^*<\infty$, $l=1,2,\ldots,M$. Hence, the minimax
rates of convergence (in the $L^2$-risk) are determined by
$\nu_{\min}$ only, meaning that one can just rely on the best
possible channel and disregard all the others. However, the latter
is no longer true if $M \rightarrow\infty$. In this case, the
minimax rates of convergence (in the $L^2$-risk) are determined by
$\tau_1 (m)$ which may not be a function of $\nu_{\min}$ only.

Consider now the adaptive block thresholding wavelet estimator
$\hat{f}_n(\cdot)$ defined by \fr{fest} for the model \fr{multchan}
$\sigma_l=1$ for $l=1,2,\ldots,M$
or its discrete counterpart (\ref{convdis}). Then, for the
$L^2$-risk, under the assumption \fr{cond1}, the corresponding
minimax lower bounds are given by Theorem \ref{th:lower}, while,
under the assumption \fr{cond2}, the corresponding minimax upper
bounds are given by Theorem \ref{th:upper}. Thus, the proposed
functional deconvolution methodology significantly expands on the
theoretical findings in De Canditiis and Pensky (\citeyear{DeCanPen2006}).
\end{example}

\section{Discussion: the interplay between continuous and discrete models}
\label{discrcont}

The minimax convergence rates (in the $L^2$-risk) in the discrete
model depend on two aspects: the total number of observations $n =
NM$ and the behavior of $\tau_1 (m)$ defined in \fr{taum}. In the
continuous model, the values of $\tau_1 (m)$ are fixed; however, in
the discrete model they may depend on the choice of $M$ and the
selection of points $u_l$, $l=1,2,\ldots, M$. Let us now explore
when and how this can happen.

Assume that there exist points $u_*, u^* \in[a,b]$, $-\infty< a
\leq b < \infty$ (with $a < b$ in the continuous model while $a=b$
is possible in the discrete model), such that $u_* = \arg\min_u
\gmu$ and $u^* = \arg\max_u \gmu$. (Obviously, this is true if the
functional Fourier coefficients $g_m(\cdot)$ are continuous
functions on the compact interval $[a,b]$.) In this case, we have
$\tau_1 (m) \geq L_* |\gm(u_*)|^2$ and $\tau_1 (m) \leq L^* |\gm
(u^*)|^2$, where $L_* = L^* = b-a$ in the continuous model and $L_*
= L^* = 1$ in the discrete model. Assume also that we can observe
$y(u,t)$ at the points $u_*$ and $u^*$. The following statement
presents the case when the minimax convergence rates cannot be
influenced by the choice of $M$ and the selection of points $u_l$,
$l=1,2,\ldots,M$.

\begin{proposition} \label{th:condisrates}
Let there exist constants $\nu_1 \in\RR$, $\nu_2 \in\RR$, $\af_1
\geq0$, $\af_2 \geq0$, $\beta_1 >0$, $\beta_2 >0$, $L_1 >0$ and
$L_2
>0$, independent of $m$, such that
%
\begin{eqnarray}
|\gm(u_*)|^2 & \geq & L_1 |m|^{-2\nu_1} \exp(-\af_1 |m|^{\beta_1}),
\qquad\nu_1 >0 \mbox{ if } \af_1 =0, \label{disconFanis}\\
|\gm(u^*)|^2 & \leq & L_2 |m|^{-2\nu_2} \exp(-\af_2 |m|^{\beta_2}),
\qquad\nu_2 >0 \mbox{ if } \af_2 =0, \label{discon}
\end{eqnarray}
where either $\af_1 \af_2 > 0$ and $\beta_1=\beta_2$ or $\af_1=
\af_2 = 0$ and $\nu_1 = \nu_2$. Then, the minimax convergence rates
obtained in Theorems \ref{th:lower} and \ref{th:upper} in the
discrete model are independent of the choice of $M$ and the
selection of points $u_l$, $l=1,2,\ldots,M$, and, hence, coincide
with the minimax convergence rates obtained in Theorems~\ref{th:lower}~and~\ref{th:upper} in the continuous model.
\end{proposition}

The validity of Proposition \ref{th:condisrates} follows trivially
from the lower and upper bounds obtained in Theorems \ref{th:lower}
and \ref{th:upper}. Proposition \ref{th:condisrates} simply states
that asymptotically (up to a constant factor) it makes absolutely no
difference whether one samples  \fr{convdis} $n$ times
at one point, say, $u_1$ or, say, $\sqrt{n}$ times at $M = \sqrt{n}$
points $u_l$. In other words, asymptotically (up to a constant
factor) each sample value $y(u_l, t_i)$, $l=1,2,\ldots, M$,
$i=1,2,\ldots, N$, gives the same amount of information and the
minimax convergence rates are not sensitive to the choice of $M$ and
the selection of points $u_l$, $l=1,2,\ldots,M$. The constants in
Theorem \ref{th:upper} will, of course, reflect the difference and
will be the smallest if one samples (\ref{convdis}) $n$ times at
$u^*$.

However, conditions \fr{disconFanis}--\fr{discon} are not always
true. Consider, for example, the case when $g(u,x) = (2u)^{-1}\,\II(|x|
\leq u)$, that is, the case of a boxcar-type convolution for each $u
\in[a,b]$, $0 <a < b < \infty$. Then, $\gmu= \sin(2 \pi m u)/(2\pi
mu)$ and $|\gm(u_*)|^2 =0$; indeed, for rational points $u=l_1/l_2
\in[a,b]$, the functional Fourier coefficients $\gm(u)$ vanish for
any integer $m$ multiple of $l_2$. This is an example where a
careful choice of $u_l$, $l=1,2,\ldots,M$, can make a difference.
For example, in the multichannel boxcar deconvolution problem (see
also Example \ref{ex5}), De Canditiis and Pensky (\citeyear{DeCanPen2006}) showed that if $M$ is
finite, $M \geq2$, one of the $u_l$'s is a ``badly approximable''
(BA) irrational number, and $u_1,u_2, \ldots, u_M$ is a BA
irrational tuple, then $\Delta_1 (j) \leq C j\, 2^{j(2+ 1/M)}$ [for
the definitions of the BA irrational number and the BA irrational
tuple, see, e.g., Schmidt (\citeyear{Sch1980})]. This implies that, in this case,
(the degree of ill-posedness is) $\nu= 1 + 1/(2M)$. [The case
$M=1$, corresponding to the standard boxcar deconvolution problem,
was considered by Johnstone, Kerkyacharian, Picard and Raimondo
(\citeyear{Joh2004}) who showed that $\nu=3/2$ when $u_1$ is a BA irrational
number.] Furthermore, De Canditiis and Pensky (\citeyear{DeCanPen2006}) obtained
asymptotical upper bounds (in the $L^{\pi}$, $1 < \pi< \infty$) for
the error, for a wavelet estimator, for a fixed response function.
They also showed that these bounds depend on $M$ and the larger the
$M$, is the higher the asymptotical convergence rates will be. Hence,
in the multichannel boxcar deconvolution problem, it is
advantageous to take $M \rightarrow\infty$ and to choose
$u_1,u_2,\ldots,u_M$ to be a BA tuple.

However, the theoretical results obtained in Theorems \ref{th:lower}
and \ref{th:upper} cannot be blindly applied to accommodate the
blurring scenario represented by the case of boxcar-type
convolution for each fixed $u$, that is, the case when $g(u,x) =
(2u)^{-1}\,\II(|x| \leq u)$, $u \in[a,b]$, $0 <a < b < \infty$. A
careful treatment of this problem is necessary, since it requires
nontrivial results in number theory. This is currently under
investigation by the authors and the results of the analysis will be
published elsewhere.

\section{Proofs}\label{append}

In what follows, for simplicity, we use the notation $g$ instead of
$g(\cdot)$, for any arbitrary function $g(\cdot)$. Also, $\psijk$
refer to the periodized Meyer wavelets defined in Section
\ref{estconstr}.

\subsection{Lower bounds}
\mbox{}
\begin{pf*}{Proof of Theorem \ref{th:lower}}
The proof of the lower bounds
falls into two parts. First, we consider the lower bounds obtained
when the worst functions $f$ (i.e., the hardest functions to
estimate) are represented by only one term in a wavelet expansion
(sparse case), and then when the worst functions $f$ are
uniformly spread over the unit interval $T$ (dense case).

\subsubsection*{Sparse case}
Consider the continuous model \fr{convcont}. Let
the functions $f_{jk}$ be of the form $\fjk= \gaj\psijk$ and let
$f_0 \equiv0$. Note that by \fr{bpqs}, in order $\fjk\in\Bpqsa$,
we need $\gaj\leq A 2^{-js'}$. Set $\gaj= c 2^{-j s'}$, where $c$
is a positive constant such that $c<A$, and apply the following
classical lemma on lower bounds:

\begin{lemma}[[H\"{a}rdle, Kerkyacharian, Picard and Tsybakov (\citeyear{HarKerPicTsy1998}), Lemma
10.1{]}] \label{korost}
Let $V$ be a functional space, and let $d(\cdot, \cdot)$
be a distance on $V$. For $f, g \in V$, denote by $\varLambda_n
(f,g)$ the likelihood ratio $\varLambda_n (f,g) = d\PP_{X_n^{(f)}}/
d\PP_{X_n^{(g)}}$, where $d\PP_{X_n^{(h)}}$ is the probability
distribution of the process $X_n$ when $h$ is true. Let $V$ contains
the functions $f_0, f_1, \ldots, f_\aleph$ such that:
\begin{longlist}
\item[\textup{(a)}] $d(f_k, f_{k'}) \geq\delta>0$ for $k=0,1,\ldots,\aleph
$,  $k \neq k'$,

\item[\textup{(b)}] $\aleph\geq\exp(\lambda_n)$ for some $\lambda_n >0$,

\item[\textup{(c)}] $\ln\varLambda_n (f_0, f_k) = u_{nk} - v_{nk}$, where
$v_{nk}$ are constants and $u_{nk}$ is a
random variable such that there exists $\pi_0>0$ with
$\PP_{f_k}(u_{nk}
>0) \geq\pi_0$,

\item[\textup{(d)}] $\sup_k v_{nk} \leq\lambda_n$.
\end{longlist}
Then, for an arbitrary estimator $\hf$,
\[
\sup_{f \in V} \PP_{X_n^{(f)}}  \bigl(d(\hf, f) \geq\delta/2 \bigr)
\geq\pi_0/2.
\]
\end{lemma}

Let now $V=\{\fjk\dvtx 0 \leq k \leq2^j-1 \}$ so that $\aleph=
2^j$. Choose $d(f,g) = \| f-g\|$, where $\| \cdot\|$ is the
$L^2$-norm on the unit interval $T$. Then, $d(\fjk, \fjkp) = \gaj=
\delta$. Let $v_{nk} = \lambda_n = j \ln2$ and $u_{nk} = \ln
\varLambda_n (f_0, \fjk) + j \ln2$. Now, to apply Lemma
\ref{korost}, we need to show that for some $\pi_0 >0$, uniformly
for all $\fjk$, we have
\[
\PP_{\fjk} (u_{nk} > 0) = \PP_{\fjk}  \bigl( \ln\varLambda_n (f_0,
\fjk) > -j \ln2 \bigr) \geq\pi_0 >0.
\]
Since, by Chebyshev's
inequality,
\[
\PP_{\fjk}  \bigl( \ln\varLambda_n (f_0, \fjk)
> -j \ln2 \bigr) \geq1 - \frac{\EE_{\fjk}  |\ln\Lam_n (f_0,
\fjk) |}{j \ln2},
\]
we need to find a uniform upper bound for
$\EE_{\fjk} |\ln\varLambda_n (f_0, \fjk)|$.

Let $W(u,t)$ and $\widetilde{W}(u,t)$ be Wiener sheets on $U \times
T$. Let $\tilde{z} (u,t) = \sqrt{n}\, (g*\fjk) (u,t) + z (u,t)$,
where $z(u,t)=\dot{W}(u,t)$ and $\tilde{z}(u,t)=\dot{\widetilde{W}}(u,t)$
[i.e., $W(u,t)$ and $\widetilde{W}(u,t)$ are the primitives of
$z(u,t)$ and $\tilde{z} (u,t)$, resp.]. 
Then, assuming that $\intT
\int_U n(g*\fjk)^2(u,t)\,du \,dt < \infty$, by the multiparameter
Girsanov formula [see, e.g., Dozzi (\citeyear{Doz1989}), page 89], we get
%
\begin{eqnarray}
\label{likelihood1}
-\ln\varLambda_n (f_0, \fjk)
&= &\sqrt{n} \int_T \int_U (g*\fjk) (u,t) \,dW(u,t) \\
&&{} - \frac{n}{2} \intT\int_U (g*\fjk)^2(u,t) \,du \,dt.\nonumber
\end{eqnarray}
Hence,
\[
\EE_{\fjk}   |\ln\varLambda_n (f_0, \fjk)  | \leq A_n +
B_n,
\]
where
\begin{eqnarray*}
A_n & = & \sqrt{n} \ga_j \EE \biggl| \intT\int_U (\psijk*g)(u,t) \,dW(u,t)  \biggr|,
\\
B_n & = & 0.5 n \ga_j^2 \intT\int_U (\psijk*g)^2 (u,t)  \,du \,dt.
\end{eqnarray*}

Since, by Jensen's inequality, $A_n \leq\sqrt{2 B_n}$, we only need
to construct an upper bound for $B_n$. For this purpose, we denote
the Fourier coefficients of $\psi(\cdot)$ by $\psi_m = \langle e_m, \psi\rangle$,
and observe that in the case of Meyer wavelets, $|\psimjk| \leq
2^{-j/2}$ [see, e.g., Johnstone, Kerkyacharian, Picard and Raimondo
(\citeyear{Joh2004}), page 565]. Therefore, by properties of the Fourier transform,
we get
%
\begin{equation}\label{bnval}
B_n = O  \Biggl( 2^{-j} n \ga_j^2 \sumcj\int_U |\gmu|^2\,du  \Biggr).
\end{equation}
Let $j=j_n$ be such that
%
\begin{equation}\frac{B_n
+ \sqrt{2 B_n}}{j \ln2} \leq\frac{1}{2}. \label{jn1}
\end{equation}
Then, by
applying Lemma \ref{korost} and Chebyshev's inequality, we obtain
\begin{eqnarray} \label{lowbound1}
\inf_{\tilde{f}_n} \sup_{f \in
\Bpqs(A)} \EE\|\tilde{f}_n - f\|^2 & \geq&
\inf_{\tilde{f}_n} \sup_{f \in V} \tfrac{1}{4} \gaj^{2}
\PP (\|\tilde{f}_n - f\| \ge\gaj/2 ) \nonumber\\[-8pt]\\[-8pt]
 & \geq&\tfrac{1}{4}
\gaj^{2} \pi_0. \nonumber
\end{eqnarray}
Thus, we just need to choose the smallest possible $j=j_n$
satisfying \fr{jn1}, to calculate $\gaj= c 2^{-j s'}$, and to plug
it into \fr{lowbound1}. By direct calculations, we derive, under
condition \fr{cond1}, that
%
\begin{equation}\label{err1}
\qquad\sumcj\int_U |\gmu|^2 \leq\cases{
C 2^{-j(2\nu-1)}, &\quad if $  \af=0$, \cr 
C 2^{-j(2\nu+ \beta-1)} \exp(-\af(2\pi/3)^\beta2^{j \beta}
), &\quad if $\af>0$, 
}
\end{equation}
so that \fr{jn1} yields $2^{j_n} = C (n/\ln n)^{1/(2 s'+ 2
\nu)}$ if $\af=0$ and $2^{j_n} = C (\ln n)^{1/\beta}$ if $\af
>0$. Hence, \fr{lowbound1} yields
%
\begin{equation}\label{lowbousparse}
\inf_{\tilde{f}_n} \sup_{f\in\Bpqs}
\EE\|\tilde{f}_n - f\|^2 \geq\cases{
C (\ln n/n)^{2 s'/(2s'+2\nu)}, &\quad if  $ \af=0$, \cr %
C (\ln n)^{-2 s'/\beta}, &\quad if $  \af>0$. 
}
\end{equation}

The proof in the discrete case is almost identical to that in the continuous
case with the only difference that [compare with \fr{likelihood1}]
\begin{eqnarray*}
- \ln\varLambda_n (f_0, \fjk) & = & 0.5 \sumN\sumM
\{[y (u_l,
t_i) - \gaj(\psijk*g)]^2(u_l, t_i) - y^2 (u_l, t_i) \}\\
& =& - v_{jk} - u_{jk},
\end{eqnarray*}
where
\begin{eqnarray*}
u_{jk} & = & \gaj \sumN\sumM(\psijk*g)(u_l, t_i) \eps_{li}, \\
v_{jk} & = & 0.5 \gaj^2 \sumN\sumM(\psijk*g)^2(u_l, t_i).
\end{eqnarray*}
Note that, due to $\PP(\eps_{li} >0) = \PP(\eps_{li} \leq0)
=0.5$, we have $\PP(u_{jk} >0) = 0.5$. Also, by properties of the
discrete Fourier transform, we get
\[
v_{jk} \leq0.5 n 2^{-j} \gaj^2  \sumcj M^{-1} \sumM|\gm(u_l)|^2.
\]
By replacing $B_n$ and $B_n + \sqrt{B_n}$ with
$v_{jk}$ in the proof for the continuous
case, and using \fr{cond1}, we arrive at \fr{lowbousparse}.

\subsubsection*{Dense case}
Consider the continuous model \fr{convcont}. Let
$\eta$ be the vector with components $\eta_k = \pm1$, $k=0,1,
\ldots, 2^j-1$, denote by $\varXi$ the set of all possible vectors
$\eta$, and let $f_{j \eta} = \gaj \sumk\eta_k \psijk$. Let also
$\eta^i$ be the vector with components $\eta^i_k =
(-1)^{\II(i=k)}\eta_k $ for $i,k = 0,1,\ldots, 2^j-1$. Note that by
\fr{bpqs}, in order $f_{j \eta} \in\Bpqsa$, we need $\gaj\leq A
2^{-j (s+1/2)}$. Set $\gaj= c_{\star} 2^{-j (s+1/2)}$, where
$c_{\star}$ is a positive constant such that $c_{\star}<A$, and
apply the following lemma on lower bounds:

\begin{lemma}[[Willer (\citeyear{Will2005}), Lemma 2{]}] \label{willer}
Let $\varLambda_n (f,g)$ be
defined as in Lemma~\ref{korost}, and let $\eta$ and $f_{j \eta}$ be
as described above. Suppose that, for some positive constants $\lambda$
and $\pi_0$, we have
\[
\PP_{f_{j \eta}} \bigl(- \ln\varLambda_n (f_{j \eta^i}, f_{j \eta})
\leq
\lambda\bigr) \geq\pi_0,
\]
uniformly for all $f_{j \eta}$ and all $i=0,\ldots, 2^j-1$. Then,
for any arbitrary estimator $\tilde{f}_n$ and for some positive constant
$C$,
\[
\max_{\eta\in\varXi} \EE_{f_{j \eta}} \|\tilde{f}_n- f_{j \eta}\| \geq C
\pi_0 e^{-\lambda}\, 2^{j/2} \ga_j.
\]
\end{lemma}

Hence, similarly to the sparse case, to obtain the lower bounds it
is sufficient to show that
\[
\EE_{f_{j \eta}} |\ln\varLambda_n
(f_{j \eta^i}, f_{j \eta})| \leq\lambda_1,
\]
for a sufficiently small
positive constant $\lambda_1$. Then, by the multiparameter Girsanov
formula [see, e.g., Dozzi (\citeyear{Doz1989}), page 89], we get
\begin{eqnarray*}
\ln\varLambda_n (f_{j \eta^i}, f_{j \eta})  &= & \sqrt{n} \int_T
\int_U \bigl(g*(f_{j \eta^i}- f_{j \eta})\bigr) (u,t) \,d W (u,t) \\
&&{} -
\frac{n}{2} \int_T \int_U \bigl(g*(f_{j \eta^i}- f_{j \eta})\bigr)^2(u,t)
\,du\,dt,
\end{eqnarray*}
and recall that $|f_{j \eta^i}- f_{j \eta})| = 2|\psi_{ji}|$. Then,
\[
\EE_{f_{j \eta}} |\ln\varLambda_n (f_{j \eta^i}, f_{j \eta})|
\leq
A_n + B_n,
\]
where
\begin{eqnarray*}
A_n &=& 2 \sqrt{n} \ga_j \EE \biggl| \intT\int_U
(\psiji*g)(u,t) \,dW(u,t)  \biggr|, \\
B_n & =& 2 n \ga_j^2 \intT\int_U (\psiji*g)^2 (u,t)  \,du \,dt.
\end{eqnarray*}
Hence, similarly to the sparse case, $A_n \leq\sqrt{2 B_n}$ and
\fr{bnval} is valid. According to Lemma \ref{willer}, we choose $j =
j_n$ that satisfies the condition $B_n + \sqrt{2 B_n} \leq\lambda_1$.
Using \fr{err1}, we derive that $2^{j_n} = C n^{1/(2s+2\nu+1)}$ if
$\af=0$ and $2^{j_n} = C (\ln n)^{1/\beta}$ if $\af>0$.
Therefore, Lemma \ref{willer} and Jensen's inequality yield
%
\begin{equation}\label{lowboudense}
\inf_{\tilde{f}_n} \sup_{f \in\Bpqs} \EE\|\tilde{f}_n - f\|^2
\geq\cases{
C n^{-2s/(2s+ 2\nu+1)}, &\quad if $ \af=0$, \cr 
C (\ln n)^{-2s/\beta}, &\quad if $  \af>0$. 
}
\end{equation}
The proof can be now extended to the discrete case in exactly
the same manner as in the sparse case. Now, to complete the proof
one just needs to note that $\sstar= \min(s,s')$, and that
%
\begin{equation}\label{cases}
2s/(2s+2\nu+1) \leq2s^*/(2s^* + 2\nu)  \qquad \mbox{if }    \nu(2-p) \leq ps^*,
\end{equation}
with the equalities taken
place simultaneously, and then to choose the highest of the lower
bounds \fr{lowbousparse} and \fr{lowboudense}. This completes the
proof of Theorem \ref{th:lower}.
\end{pf*}

\subsection{Upper bounds}

\mbox{}
\begin{pf*}{Proof of Lemma \ref{l:coef}}
In what follows, we shall only
construct the proof for $\bjk$ [i.e., the proof of \fr{hb}] since
the proof for $\ajk$ [i.e., the proof of \fr{ha}] is very similar.
First, consider the continuous model \fr{convcont}. Note that, by
\fr{coefest},
\[
\hbjk- \bjk= \sum_{m \in C_j} (\hfm- \fm) \overline{\psimjk},
\]
where
%
\begin{equation}\label{conterror}
\hfm- \fm= n^{-1/2} \biggl(\iab\bgmu \zmu \,du \biggr) \Big/ \biggl(
\iab|\gmu|^2 \,du \biggr),
\end{equation}
due to \fr{finaleq} and
\fr{fmexprc}. Recall that $\zmu$ are Gaussian processes with zero
mean and covariance function satisfying \fr{zmu}. Hence, it is easy
to check that
\[
\EE [(\widehat{f}_{m_1} - f_{m_1})\overline{(\widehat{f}_{m_2}
- f_{m_2})} ] =
n^{-1} [\tau_1 (m_1)]^{-1} \delta(m_1 -m_2),
\]
implying that
\[
\EE|\hbjk- \bjk|^2 = n^{-1} \sum_{m \in C_j} |\psimjk|^2 [\tau_1
(m)]^{-1},
\]
where $\tau_1 (m)$ is defined in (\ref{taum}) (the
continuous case). To complete the proof of \fr{hb} in the case of
$\kappa=1$, just recall that $|C_j| = 4 \pi2^j $ and $|\psimjk|^2
\leq2^{-j}$. If $\kappa=2$, then
\begin{eqnarray*}
\EE|\hbjk-\bjk|^4 &= & O  \Biggl( \sum_{m \in C_j} \EE|\hfm- \fm|^4  \Biggr)
+
O  \Biggl(  \Biggl[ \sum_{m \in C_j} \EE|\hfm- \fm|^2  \Biggr]^2
\Biggr) \\
&= & O  \Biggl( n^{-2} \sum_{m \in C_j} |\psimjk|^4  \tau_2 (m) [\tau_1 (m)]^{-4}  \Biggr) \\
&&{} + O  \Biggl( n^{-2}  \Biggl[ |C_j|^{-1} \sum_{m \in C_j} [\tau_1
(m)]^{-1}  \Biggr]^2  \Biggr)\\
&= & O (n^{-2} 2^{-j} \Delta_2 (j) )+ O (n^{-2}
\Delta_1^2 (j) )= O (n^{-2} \Delta_2 (j) ),
\end{eqnarray*}
since, by the Cauchy--Schwarz inequality, $\Delta_1^2
(j) \leq\Delta_2 (j)$. This completes the proof of \fr{hb} in the
continuous case.

In the discrete case, formula \fr{conterror} takes the form [see
\fr{fmexprd}]
%
\begin{equation}\label{discerror}
\hfm- \fm= N^{-1/2}  \Biggl( \sumM
\overline{\gm(u_l)}\, z_{ml}  \Biggr)  \bigg/  \Biggl( \sumM|\gm(u_l)|^2
 \Biggr),
\end{equation}
where $z_{ml}$ are standard Gaussian
random variables, independent for different $m$ and $l$. Therefore,
similarly to the continuous case,
\[
\EE|\hbjk-
\bjk|^2 = N^{-1} \sum_{m \in C_j} |\psimjk|^2  \Biggl[ \sumM
|\gm(u_l)|^2  \Biggr]^{-1} = O  ( n^{-1} \Delta_1 (j)  ).
\]
In the case of $\kappa=2$, note that
\[
\EE|\hbjk- \bjk|^4 = O  \bigl( 2^{-j} N^{-2} M^{-3} \Delta_2 (j) +
N^{-2} M^{-2} \Delta_1^2 (j)  \bigr) = O  ( n^{-2} \Delta_2 (j) ),
\]
by applying again the Cauchy--Schwarz inequality. This
completes the proof of \fr{hb} in the discrete case.

The last part of the lemma follows easily from (\ref{Delt}) with
$\kappa=2$, using the assumption (\ref{cond1}) and the
Cauchy--Schwarz inequality, thus completing the proof of Lemma \ref{l:coef}.
\end{pf*}

\begin{pf*}{Proof of Lemma \ref{l:deviation}}
Consider the set of vectors
\[
\varOmega_{jr} =  \Biggl\{ v_k,  k \in\Ujr\dvtx \sumku|v_k|^2 \leq1\Biggr\},
\]
and the centered Gaussian process defined by
\[
Z_{jr} (v) = \sumku v_k (\hbjk- \bjk).
\]
The proof of the lemma is based on the following inequality:

\begin{lemma}[[Cirelson, Ibragimov and Sudakov (\citeyear{CirIbrSud1976}){]}]\label{l:cirel}
Let $D$ be a subset
of $\RR=(-\infty,\infty)$, and let $(\xi_t)_{t \in D}$ be a
centered Gaussian process. If~ $\EE(\sup_{t \in D} \xi_t ) \leq B_1$
and $\sup_{t \in D}\operatorname{Var} (\xi_t) \leq B_2$, then, for all $x>0$,
we have
%
\begin{equation}\label{cirel}
\PP \biggl( \sup_{t \in D} \xi_t \geq x + B_1 \biggr)
\leq
\exp (-x^2/(2 B_2) ).
\end{equation}
\end{lemma}

To apply Lemma~ \ref{l:cirel}, we need to find $B_1$ and $B_2$. Note
that, by Jensen's inequality, we get
\begin{eqnarray*}
\EE \biggl[\sup_{v \in\varOmega_{jr}} Z_{jr} (v)  \biggr]
& = & \EE \Biggl[ \sumku
|\hbjk- \bjk|^2
 \Biggr]^{1/2}\\
  & \leq & \Biggl[ \sumku\EE|\hbjk- \bjk|^2  \Biggr]^{1/2}\\
& \leq&\sqrt{c_1} n^{-1/2} 2^{\nu j} \sqrt{\ln n}.
\end{eqnarray*}
[Here, $c_1$ is the same positive constant as in \fr{delta1} with
$\alpha=0$.] Also, by \fr{zmu} and \fr{conterror} or \fr{discerror},
we have
\[
\EE [(\hbjk- \bjk)(\widehat{b}_{jk'} - b_{jk'})  ] = n^{-1}
\sumcj\psimjk
\overline{\psi_{mjk'}} [\tau_1 (m)]^{-1},
\]
where $\tau_1 (m)$ is
defined in (\ref{taum}). Hence,
\begin{eqnarray*}
\sup_{v \in\varOmega_{jr}} \Var(Z_{jr} (v)) & = & n^{-1} \sup_{v \in
\varOmega_{jr}} \sumku\sum_{k' \in\Ujr} v_k v_{k'}
\sumcj\psimjk\overline{\psi_{mjk'}} [\tau_1 (m)]^{-1} \\
& \leq&  c_1 n^{-1} 2^{2 \nu j} \sumku v^2_k \leq c_1 n^{-1} 2^{2 \nu
j},
\end{eqnarray*}
by using $\sumcj\psimjk\overline{\psi_{mjk'}} =
\II(k=k')$ and \fr{delta1} for $\alpha=0$. Therefore, by applying
Lemma~ \ref{l:cirel} with $B_1 = \sqrt{c_1} n^{-1/2} 2^{\nu j}
\sqrt{\ln n}$, $B_2 = c_1 n^{-1} 2^{2 \nu j}$ and $x = (0.5 \mu-
\sqrt{c_1}) n^{-1/2} 2^{\nu j} \sqrt{\ln n}$, we get
\begin{eqnarray*}
&&\PP \Biggl( \sumku|\hbjk- \bjk|^2    \geq0.25
\mu^2 n^{-1} 2^{2 \nu j} \ln n  \Biggr) \\
&&\qquad =
\PP \Biggl(  \Biggl[ \sumku|\hbjk- \bjk|^2  \Biggr]^{1/2} \geq\sqrt
{c_1} n^{-1/2} 2^{\nu j} \sqrt{\ln n} + x  \Biggr) \\
&&\qquad \leq\exp\bigl(- (2 c_1)^{-1} \bigl(0.5 \mu- \sqrt{c_1}\bigr)^2  \ln
n\bigr) \leq n^{-\te},
\end{eqnarray*}
where $\te= (8\nu- 4 \nu_1
+2)/(2\nu+1)$, provided that $\mu\geq2 \sqrt{c_1} (1 + \sqrt{2
\te})$. This completes the proof of Lemma \ref{l:deviation}.
\end{pf*}

\begin{pf*}{Proof of Theorem \ref{th:upper}}
First, note that in the case
of $\af>0$, we have
\[
\EE\| \hf- f \|^2 = R_1 + R_2,
\]
where
%
\begin{equation}\label{r1r2}
R_1 = \sum_{j=J}^\infty\sumk\bjk^2, \qquad   R_2 =
\sum_{k=0}^{2^{j_0}-1} \EE(\hajk- \ajk)^2,
\end{equation}
since
$\jo= J$. It is well known [see, e.g., Johnstone (\citeyear{Joh2002}), Lemma
19.1] that if $f \in\Bpqsa$, then for some positive constant
$c^{\star}$, dependent on $p$, $q$, $s$ and $A$ only, we have
%
\begin{equation}
\sumk\bjk^2 \leq c^{\star} 2^{-2 j \sstar}, \label{besball}
\end{equation}
thus, $R_1 = O (2^{-2 J \sstar} )= O ((\ln
n)^{-2\sstar/\beta} )$. Also, using \fr{delta1} and \fr{ha}, we
derive
\[
R_2 = O (n^{-1} 2^{j_0}\Delta_1 (\jo) )= O (n^{-1/2}
(\ln n)^{2\nu/\beta} )= o ((\ln n)^{-2\sstar/\beta} ),
\]
thus completing the proof for $\af>0$.

Now, consider the case of $\af=0$. Due to the orthonormality of
the wavelet basis, we get
%
\begin{equation}\label{errtotal}
\EE\| \hf- f \|^2 = R_1 + R_2 + R_3+ R_4,
\end{equation}
where $R_1$ and $R_2$ are defined in
\fr{r1r2}, and
\begin{eqnarray*}
R_3 & = & \sumjJ\sumr\sumku\EE[(\hbjk- \bjk)^2 \II(\hBjr
\geq d n^{-1} 2^{2\nu j} \ln n ) ], \\
R_4 & = & \sumjJ\sumr\sumku\EE[\bjk^2 \II(\hBjr< d n^{-1}
2^{2\nu j} \ln n ) ],
\end{eqnarray*}
where $\hBjr$ and
$d$ are given by \fr{bjr} and \fr{lamj}, respectively.

Let us now examine each term in \fr{errtotal} separately. Similarly
to the case of $\af>0$, we obtain $R_1 = O (2^{-2 J \sstar}
)= O (n^{-2\sstar/(2 \nu+1)} ).$ By direct calculations,
one can check that $2\sstar/(2 \nu+1) > 2s/(2s + 2\nu+1)$, if
$\nu(2-p) < p \sstar$, and $2\sstar/(2 \nu+1) \ge2\sstar/(2\sstar
+ 2\nu)$, if $\nu(2-p) \geq p \sstar$. Hence,
%
\begin{eqnarray}
R_1 &=& O \bigl(n^{-2s/(2s+2\nu+1)} \bigr)\qquad  \mbox{if }   \nu
(2-p) < p \sstar, \label{r1a}\\
R_1 &=& O \bigl(n^{-2\sstar/(2\sstar+2\nu)} \bigr)\hspace*{2.6pt} \qquad\mbox{if }   \nu(2-p) \geq p \sstar. \label{r1b}
\end{eqnarray}
Also, by \fr{ha} and \fr{delta1}, we get
\begin{eqnarray}\label{r2}
R_2 & = & O \bigl(n^{-1}
2^{(2 \nu+1) \jo} \bigr)= O (n^{-1} (\ln n)^{2\nu+1}
)\nonumber\\[-8pt]\\[-8pt]
& = & o \bigl(n^{-2s/(2s+2\nu+1)} \bigr)= o \bigl(
n^{-2s^*/(2s^*+2\nu)} \bigr). \nonumber
\end{eqnarray}

To construct the upper bounds for $R_3$ and $R_4$, note that simple
algebra gets
%
\begin{equation}
R_3 \leq(R_{31} + R_{32}), \qquad   R_4 \leq(R_{41}
+ R_{42}), \label{r34}
\end{equation}
where
\begin{eqnarray*}
R_{31} & = & \sumjJ
\sumr\sumku\EE\Biggl[(\hbjk- \bjk)^2 \II \Biggl( \sumku|\hbjk-
\bjk|^2 \geq0.25 d n^{-1} 2^{2 \nu j}
\ln n  \Biggr) \Biggr], \nonumber\\
R_{32} & = & \sumjJ\sumr\sumku\EE[(\hbjk- \bjk)^2  \II
(\Bjr> 0.25 d n^{-1} 2^{2 \nu j} \ln n  ) ], \nonumber\\
R_{41} & = & \sumjJ\sumr\sumku\EE\Biggl[\bjk^2 \II \Biggl( \sumku
|\hbjk- \bjk|^2 \geq0.25 d n^{-1} 2^{2 \nu j} \ln n  \Biggr) \Biggr],
\nonumber\\
R_{42} & = & \sumjJ\sumr\sumku\EE[\bjk^2 \II (\Bjr< 2.5
d n^{-1} 2^{2 \nu j} \ln n  ) ], \nonumber
\end{eqnarray*}
since $\hbjk^2 \leq2(\hbjk- \bjk)^2 + 2\bjk^2$. Then, by
\fr{besball}, Lemmas \ref{l:coef} and \ref{l:deviation}, and the
Cauchy--Schwarz inequality, we derive
\begin{eqnarray*}
R_{31}+ R_{41}
& =& \sumjJ\sumr\sumku\EE \Biggl( \bigl((\hbjk- \bjk)^2 + \bjk^2
\bigr)\\
&&\hspace*{79pt}
{}\times \II \Biggl( \sumku|\hbjk- \bjk|^2
\geq0.25 d n 2^{2 \nu j} \ln n  \Biggr)  \Biggr) \\
& = & O  \Biggl( \sumjJ\sumr\sumku \bigl( \sqrt{\EE(\hbjk-
\bjk)^4} + \bjk^2  \bigr)\\
&&\phantom{O  \Biggl(}\hspace*{63pt}
{}\times
\sqrt{ \PP \Biggl( \sumku|\hbjk- \bjk|^2 \geq0.25 d n 2^{2 \nu j}
\ln n  \Biggr)}  \Biggr) \\
& = & O \Biggl(\sumjJ\bigl[2^j n^{-1} 2^{2(2\nu-\nu_1)j} + 2^{-2j\sstar}\bigr] n^{-(4\nu- 2 \nu_1 + 1)/(2\nu+1)} \Biggr)\\
&=& O(n^{-1}),
\end{eqnarray*}
provided $d \geq c_1 (1 + \sqrt{2 \te})^2$, where $\te= (8\nu- 4
\nu_1 +2)/(2\nu+1)$ and $c_1$ is the same positive constant as in
\fr{delta1} with $\alpha=0$. Hence,
%
\begin{equation}\varDelta_1 = R_{31}+ R_{41}
= O (n^{-1} ). \label{r3141}
\end{equation}

Now, consider
%
\begin{equation}
\varDelta_2= R_{32} + R_{42}. \label{del2}
\end{equation}
Let $j_1$ be such that
%
\begin{equation}
2^{j_1} = n^{1/(2s + 2\nu+1)} (\ln
n)^{\ro_1/(2\nu+1)}, \label{j1}
\end{equation}
where $\ro_1$ is defined
in \fr{rovalue}.

First, let us study the dense case, that is, when $\nu(2-p) < p{s^*}$.
Then, $\varDelta_2$ can be partitioned as $\varDelta_2 =
\varDelta_{21} + \varDelta_{22}$, where the first component is
calculated over the set of indices $\jo\leq j \leq j_1$ and the
second component over $j_1+1 \leq j \leq J-1$. Hence, using \fr{bjr}
and Lemma \ref{l:coef}, and taking into account that the cardinality
of $A_j$ is $|A_j|=2^j/\ln n$, we obtain
\begin{eqnarray}
\label{del21c1}
\varDelta_{21} & =& O  \Biggl( \sumjone \Biggl[ \sumk\EE(\hbjk- \bjk
)^2 \nonumber\\
& &\phantom{O  \Biggl( \sumjone \Biggl[}
{}+ \sumr\Bjr\II (\Bjr\leq2.5 d n^{-1} 2^{2 \nu j}
\ln n  )  \Biggr]  \Biggr) \nonumber\\[-8pt]\\[-8pt]
& = & O \Biggl(\sumjone \Biggl[ n^{-1} 2^{(2\nu+1)j} + \sumr n^{-1} 2^{2
\nu j} \ln n  \Biggr] \Biggr)\nonumber\\
& = & O \bigl(n^{-1} 2^{ (2\nu+1)j_1} \bigr)= O \bigl(n^{-2s/(2s
+ 2\nu+1)} (\ln n)^{\ro_1} \bigr).\nonumber
\end{eqnarray}
To obtain an expression for $\varDelta_{22}$, note that, by
\fr{besball}, and for $p \geq2$, we have
%
\begin{eqnarray}
\label{del22c0}
\varDelta_{22} & =& O  \Biggl( \sumjtwo\sumr [ n^{-1} 2^{2\nu
j} \ln n    \mathbb{I}  ( \Bjr
\geq0.25 d n^{-1} 2^{2 \nu j} \ln n  )
+ \Bjr ]  \Biggr)\nonumber \\
& = & O  \Biggl( \sumjtwo\sumr\Bjr \Biggr) = O  \Biggl( \sumjtwo
2^{-2js}  \Biggr)
\\
& = & O  \bigl( n^{-2s/(2s+2\nu+1)}  \bigr).\nonumber
\end{eqnarray}
If $1 \leq p <2$, then
\[
\Bjr^{p/2} = \Biggl(\sumku\bjk^2 \Biggr)^{p/2} \leq\sumku|\bjk|^p,
\]
so that by Lemma \ref{l:coef}, and since $\nu(2 - p) < p s^*$, we
obtain
%
\begin{eqnarray}
\label{del22c1}
\varDelta_{22} & = & O  \Biggl( \sumjtwo\sumr [ n^{-1} 2^{2\nu
j}
\ln n  \II (\Bjr\geq0.25 d n^{-1} 2^{2 \nu j} \ln n  )\nonumber \\
 & &\hspace*{97pt}{}+ \Bjr\II (\Bjr\leq2.5 d n^{-1} 2^{2 \nu j} \ln
n )  ]  \Biggr)
\\
& = & O  \Biggl( \sumjtwo\sumr [  ( n^{-1} 2^{2 \nu j} \ln n
 )^{1-p/2} \Bjr^{p/2} \nonumber\\
 &&\hspace*{69pt}{} + \Bjr^{p/2}  ( n^{-1} 2^{2
\nu j} \ln n  )^{1-p/2} ]  \Biggr)\nonumber
\\
& =& O \Biggl(\sumjtwo ( n^{-1} 2^{2 \nu j} \ln n
)^{1-p/2} \sumr\sumku|\bjk|^p \Biggr)\nonumber
\\
& =& O \Biggl(\sumjtwo ( n^{-1} 2^{2 \nu j} \ln n
)^{1-p/2} 2^{-pjs^*} \Biggr)\nonumber
\\
& =& O \Biggl(\sumjtwo(n^{-1} \ln n )^{1 - p/2} 2^{(2 \nu- p\nu-
p s^*)j} \Biggr)\nonumber
\\
& =& O \bigl((n^{-1} \ln n )^{1 - p/2}  2^{ (2 \nu- p\nu- p
s^*)j_1} \bigr)\nonumber
\\
& =& O \bigl(n^{-2s/(2s+2\nu+1)} (\ln n)^{\ro_1} \bigr).\nonumber
\end{eqnarray}

Let us now study the sparse case, that is, when $\nu(2-p) > p{s^*}$.
Let $j_1$ be defined by \fr{j1} with $\ro_1 =0$. Hence, if $\Bjr
\geq0.25 d n^{-1} 2^{2 \nu j} \ln n $, then $\sumk\bjk^2 \geq0.25
d n^{-1} 2^{(2 \nu+1)j}$, implying that $j$ cannot exceed $j_2$
such that $2^{j_2} = (4c^*n/(d\,{\ln n}))^{1/(2\sstar+ 2\nu)}$,
where $c^*$ is the same constant as in \fr{besball}. Again,
partition $\varDelta_2 = \varDelta_{21} + \varDelta_{22}$, where the
first component is calculated over $\jo\leq j \leq j_2$ and the
second component over $j_2+1 \leq j \leq J-1$. Then, using arguments similar to those in \fr{del22c1}, and taking into account that $
\nu(2 - p)>p \sstar$, we derive
%
\begin{eqnarray}
\label{del21c2}
\varDelta_{21} & = & O \Biggl(\sumjones ( n^{-1}2^{2 \nu j} \ln n
 )^{1-p/2} \sumr\sumku|\bjk|^p \Biggr)\nonumber\\
 & = & O \Biggl(\sumjones ( n^{-1} 2^{2 \nu j} \ln n  )^{1-p/2}
2^{-pj\sstar} \Biggr)\nonumber\\
& = & O \Biggl(\sumjones(n^{-1} \ln n)^{1 - p/2} 2^{(2 \nu- p\nu-
p \sstar)j} \Biggr)\\
& = & O \bigl((n^{-1} \ln n )^{1 - p/2}  2^{(2 \nu- p\nu- p \sstar
)j_2 } \bigr)\nonumber\\
& = & O \bigl( (\ln n/n )^{2{s^*}/(2s^*+2\nu)} \bigr).\nonumber
\end{eqnarray}
To obtain an upper bound for $\varDelta_{22}$, recall \fr{del2} and
keep in mind that the portion of $R_{32}$ corresponding to $j_2+1
\leq j \leq J-1$ is just zero. Hence, by \fr{besball}, we get
\begin{eqnarray} \label{del22c2}
\varDelta_{22} &=& O \Biggl(\sumjtwos\sumk
\bjk^2 \Biggr) = O \Biggl(\sumjtwos2^{-2j\sstar} \Biggr)
\nonumber\\[-8pt]\\[-8pt]
& =& O \bigl( (\ln n/n  )^{2{s^*}/(2s^*+2\nu)} \bigr).
\nonumber
\end{eqnarray}

Now, in order to complete the proof, we just need to study the case
when $\nu(2-p) = p{s^*}$. In this situation, we have $2s/(2s+ 2\nu
+ 1) = 2s^* /(2s^*+ 2\nu) = 1-p/2$ and $2 \nu j (1 - p/2) = p j
s^*$. Recalling \fr{bpqs} and noting that $s^* \leq s'$, we get
\[
\sumjJ\Biggl(2^{p j s^*} \sumk|\bjk|^p \Biggr)^{q/p} \leq A^q.
\]
Then, we repeat the calculations in \fr{del21c2} for all indices
$\jo\leq j \leq J-1$. If $1 \leq p < q$, then, by H\"{o}lder's
inequality, we get
%
\begin{eqnarray}
\label{del2c3}
\varDelta_{2} & =& O  \Biggl( \sumjJ(\ln n/n )^{1-p/2} 2^{p j
s^*} \sumk|\bjk|^p  \Biggr) \nonumber\\
& = & O  \Biggl(  ( \ln n/n  )^{1-p/2} (\ln n)^{1 - p/q}  \Biggl[ \sumjJ \Biggl(
2^{p j s^*} \sumk|\bjk|^p  \Biggr)^{q/p}  \Biggr]^{p/q}  \Biggr) \\
& = & O  \bigl(  (\ln n/n  )^{2{s^*}/(2s^*+2\nu)} (\ln
n)^{1 - p/q}  \bigr).\nonumber
\end{eqnarray}
If $1 \leq q \leq p$, then, by the inclusion $\Bpqs(A) \subset
B_{p,p}^s(A)$, we get
%
\begin{eqnarray}
\label{del2c4FF}
\varDelta_{2} & = & O  \Biggl( \sumjJ(\ln n/n )^{1-p/2} 2^{p j
s^*} \sumk|\bjk|^p  \Biggr)\nonumber \\
& = & O  \Biggl(  ( \ln n/n  )^{1-p/2}   \Biggl[ \sumjJ2^{p j
s^*} \sumk|\bjk|^p  \Biggr]  \Biggr) \\
& = & O  \bigl(  (\ln n/n  )^{2{s^*}/(2s^*+2\nu)}\bigr).\nonumber
\end{eqnarray}
By combining \fr{r1a}--\fr{r2}, \fr{r3141},
\fr{del21c1}--\fr{del2c4FF}, we complete the proof of Theorem
\ref{th:upper}.
\end{pf*}

\section*{Acknowledgments}

Marianna Pensky is grateful for the hospitality and financial
support of the Department of Mathematics and Statistics at the
University of Cyprus, Cyprus, and Theofanis Sapatinas is grateful
for the hospitality of the Department of Mathematics at the
University of Central Florida, USA, where parts of the work of this
paper were carried out. The authors would like to thank Thomas
Willer for useful discussions. Finally, we would like to thank an
Associate Editor and two anonymous referees for their suggestions on
improvements to this paper.

\printaddresses

\end{document}